\documentclass[pdflatex,sn-mathphys-num]{sn-jnl}

\usepackage{caption}
\usepackage{subcaption}
\usepackage{bm}

\usepackage{enumitem}

\usepackage{graphicx}
\usepackage{float}

\usepackage{booktabs}
\usepackage{tablefootnote}

\usepackage{mathtools}
\usepackage{csquotes}

\usepackage{amsmath}
\usepackage{amssymb}
\usepackage{amsfonts}

\usepackage{cleveref}

\usepackage{multirow}

\usepackage{tikz}
\usepackage{pgfplots}
\usetikzlibrary{patterns}
\pgfplotsset{compat=1.18}
\usetikzlibrary{trees,shapes,decorations}
\usepgfplotslibrary{groupplots}

\input{macros}

\usepackage{lipsum}
\usepackage{epstopdf}
\ifpdf
  \DeclareGraphicsExtensions{.eps,.pdf,.png,.jpg}
\else
  \DeclareGraphicsExtensions{.eps}
\fi

\usepackage[ruled,vlined,algo2e]{algorithm2e} 

\SetKwProg{Function}{function}{}{}
\SetKw{KwRet}{return}
\DontPrintSemicolon
\makeatother

\theoremstyle{thmstyleone}
\newtheorem{theorem}{Theorem}
\newtheorem{lemma}[theorem]{Lemma}
\newtheorem{corollary}[theorem]{Corollary}
\newtheorem{proposition}[theorem]{Proposition}

\theoremstyle{thmstyletwo}%
\newtheorem{example}{Example}%
\newtheorem{remark}{Remark}%

\theoremstyle{thmstylethree}%
\newtheorem{definition}{Definition}%

\begin{document}

\title[Hierarchical splitting methods]{A hierarchical splitting approach for N-split differential equations}

\author*[1]{\fnm{Kevin} \sur{Schäfers}}\email{schaefers@math.uni-wuppertal.de}

\author[1]{\fnm{Michael} \sur{Günther}}\email{guenther@uni-wuppertal.de}

\affil*[1]{\orgdiv{IMACM}, \orgname{Bergische Universität Wuppertal}, \orgaddress{\street{Gaußstr.~20}, \postcode{42119} \city{Wuppertal}, \country{Germany}}}

\abstract{We propose a hierarchical splitting approach that constructs splitting methods for $\Npartitions$-split differential equations by iteratively applying splitting methods for two-split systems along a splitting tree.
Since assigning methods of order $\convergenceorder$ to all nodes of the tree yields a hierarchical splitting method of order $\convergenceorder$, the order theory for two-split systems carries over to any number of partitions.
We derive explicit formulas for the leading-order error terms, investigate self-adjointness, and show that state-of-the-art composition approaches are contained in the resulting class as special cases.
Augmenting the approach with multiple time stepping yields a framework at the intersection of geometric numerical integration and multirate integration, for which we characterize the computational order of a multirate integrator and introduce a reweighting of the multirate factors that reduces computational overhead, but can lead to an order reduction when applied within composition techniques.
We further discuss principles for designing splitting trees, with particular attention to separable Hamiltonian systems and port-Hamiltonian systems.
Numerical experiments for both classes demonstrate that the hierarchical splitting approach yields efficient methods for $\Npartitions$-split systems.
We find that the approach is particularly efficient when commutator-based splitting methods become available and the partitions evolve on distinct time scales.
}



\maketitle

\section{Introduction}\label{sec:introduction}
In many technical applications, such as in multi-physics problems, initial value problems (IVPs) of \emph{$\Npartitions$-split} ordinary differential equations (ODEs) of the form
\begin{equation}\label{eq:split-system}
\begin{aligned}
    \x^{\prime} &= \rhs{}(\x) = \sum_{m=1}^\Npartitions \rhs{m}(\x), & \x(0) &= \x_0,
\end{aligned}
\end{equation}
arise.
In \eqref{eq:split-system}, the right-hand side is split additively into $\Npartitions \geq 2$ partitions, based on characteristics such as stiffness, nonlinearity, dimension, dynamical behavior, and computational cost.
In such cases, it is advantageous to employ splitting approaches to numerically integrate the individual partitions by utilizing dedicated methods.
For instance, a splitting into stiff and non-stiff parts enables the numerical treatment with implicit-explicit (IMEX) methods \cite{ASCHER1997151}.
Several techniques already exist that allow the simultaneous use of dedicated numerical integration schemes for~\eqref{eq:split-system}.
For instance, in the Runge--Kutta (RK) framework, there are partitioned Runge--Kutta (PRK) methods \cite{Hairer1981PRK}, additive RK (ARK) methods \cite{cooper1983additive,KENNEDY2003139}, and generalized ARK (GARK) methods \cite{sandu2015gark}.
Another class of methods is given by splitting methods~\cite{McLachlan_Quispel_2002,blanes2024splitting}.
In contrast to RK-based methods, splitting methods are not limited to RK methods for numerical integration of the subsystems.
On the other hand, the RK-based frameworks also allow for coupling subsystems for a sub-step which may result in superior stability properties, as it is done for example in the compound step approach \cite{el2006aspects} in multirate time integration.
If the flows of all subsystems are approximated via RK methods, a splitting method can be formulated as a GARK method.
This subclass of GARK methods is also known as fractional-step RK (FSRK) methods~\cite{SPITERI2023111900}.

In addition to numerous works on splitting methods for two-split systems, there are also some works that deal with splitting methods for $\Npartitions > 2$ partitions, see e.g.~\cite{spiteri_beyond_2025,SPITERI2026115271}.
Given a step size $\stepsize > 0$, the Lie--Trotter splitting \cite{trotter_product_1959} and the Strang splitting \cite{strang1968construction}
\begin{align}\label{eq:Lie-Trotter_N2}
    \SplittingMethod{\stepsize}{} &=  \flow{\stepsize}{2} \circ \flow{\stepsize}{1}, \quad \text{and} \\
    \SplittingMethod{\stepsize}{} &= \flow{\stepsize/2}{1} \circ \flow{\stepsize}{2} \circ \flow{\stepsize/2}{1}, \label{eq:Strang_N2}
\end{align}
can be generalized to $\Npartitions \geq 2$ partitions, resulting in the methods
\begin{align}\label{eq:Lie-Trotter}
    \SplittingMethod{\stepsize}{} &= \flow{\stepsize}{\Npartitions} \circ \ldots \circ \flow{\stepsize}{2} \circ \flow{\stepsize}{1}, \quad \text{and} \\
    \SplittingMethod{\stepsize}{}&=
    \flow{\stepsize/2}{1} \circ \flow{\stepsize/2}{2} \circ \ldots \circ \flow{\stepsize/2}{\Npartitions-1} \circ \flow{\stepsize}{\Npartitions} \circ \flow{\stepsize/2}{\Npartitions-1} \circ \ldots \circ \flow{\stepsize/2}{2} \circ \flow{\stepsize/2}{1}. \label{eq:Strang-splitting}
\end{align}
Based on the methods \eqref{eq:Lie-Trotter} and \eqref{eq:Strang-splitting}, we can derive splitting methods for $\Npartitions$-split systems of arbitrarily high convergence order using composition techniques.
In one approach, we start with a base scheme $\SplittingMethod{\stepsize}{}$ of even order (such as the Strang splitting \eqref{eq:Strang-splitting}), coefficients $\gamma_1,\ldots,\gamma_s \in \mathbb{R}$, and consider the composition
\begin{equation}\label{eq:composition_based_on_Strang}
    \Psi_\stepsize = \SplittingMethod{\gamma_s \stepsize}{} \circ \ldots \circ \SplittingMethod{\gamma_1 \stepsize}{}.
\end{equation}
Alternatively, we start from a first-order base method $\chi_{\stepsize}$ (such as the Lie--Trotter splitting \eqref{eq:Lie-Trotter}) and its adjoint $\chi_{\stepsize}^*$, giving rise to composition schemes of the form
\begin{equation}\label{eq:composition_based_on_LT}
    \Psi_\stepsize = \chi_{\beta_s \stepsize} \circ \chi_{\alpha_s \stepsize}^{*} \circ \ldots \circ \chi_{\beta_1 \stepsize} \circ \chi_{\alpha_1 \stepsize}^{*}.
\end{equation}
Methods of the form \eqref{eq:composition_based_on_Strang} and \eqref{eq:composition_based_on_LT} designed for $\Npartitions = 2$ partitions and for general $\Npartitions$-split problems have the same order \cite{mclachlan1995}, but the optimization of the integrator coefficients may be different.
Unfortunately, composition schemes often have large truncation errors~\cite{omelyan2003symplectic}.
It is therefore preferable to derive splitting methods by directly solving the order conditions, which, as per Gröbner’s Lemma \cite[III Lemma 5.1]{HairerLubichWanner}, can be derived from the Baker--Campbell--Hausdorff formula for $\Npartitions = 2$ partitions or via its recent generalization to $\Npartitions \geq 2$ partitions \cite{SPITERI2026115271}.
For the GARK framework, a comprehensive order theory derived from NB-series \cite{araujo1997decomposition} is also available \cite{sandu2015gark}.
Both approaches share a common limitation: for $\Npartitions$-split systems, the number of order conditions increases significantly with $\Npartitions$, turning the construction of higher-order methods into a computationally challenging task.
Except for the compositions \eqref{eq:composition_based_on_Strang} and \eqref{eq:composition_based_on_LT} based on the Strang splitting \eqref{eq:Strang-splitting} and the Lie--Trotter splitting \eqref{eq:Lie-Trotter}, there is currently no well-established approach for designing splitting methods that can be adapted to a varying number of partitions $\Npartitions$.

This paper presents a novel class of splitting methods for $\Npartitions$-split systems referred to as \emph{hierarchical splitting methods}.
These methods apply splitting methods for two-split systems in a hierarchical manner.
Since assigning methods of order $\convergenceorder$ to all nodes of the underlying splitting tree yields a hierarchical splitting method of order $\convergenceorder$, the order theory for two-split systems extends to any number of partitions $\Npartitions$.
Furthermore, we will include multiple time stepping techniques \cite{tuckerman_reversible_1992,SEXTON1992665,URBACH200687}, placing the approach at the intersection of geometric numerical integration and multirate integration \cite{schaefers2023symplecticMGARK,oberblöbaum2024variational}.
We also show that splitting methods for $\Npartitions$-split systems that are based on the composition approaches \eqref{eq:composition_based_on_Strang} and \eqref{eq:composition_based_on_LT} can be regarded as special cases of the resulting class, and discuss design principles for splitting trees, with particular attention to nested integration schemes for separable Hamiltonian systems and to splitting approaches for port-Hamiltonian systems.

The remainder of the paper is structured as follows.
\Cref{sec:HSM} formally defines hierarchical splitting methods, discusses their convergence order, and provides explicit formulas for the principal error terms.
\Cref{sec:multiple_time_stepping} extends the approach by multiple time stepping techniques and characterizes the computational order of a multirate integrator.
\Cref{sec:design} is devoted to the design of splitting trees, and \Cref{sec:NumericalResults} presents numerical experiments for a Fermi--Pasta--Ulam problem and for coupled nonlinear mass-spring-damper chains. 

\section{Hierarchical splitting methods}\label{sec:HSM}
In this section, we will formally define the class of hierarchical splitting methods that constitutes a subclass of general $s$-stage splitting methods 
\begin{equation}\label{eq:splitting_method}
    \SplittingMethod{\stepsize}{} = \flow{\alpha_s^{\{\Npartitions\}}\stepsize}{\Npartitions} \circ \flow{\alpha_s^{\{\Npartitions-1\}}\stepsize}{\Npartitions-1} \circ \ldots \circ \flow{\alpha_s^{\{1\}}\stepsize}{1} \circ \ldots \circ \flow{\alpha_1^{\{\Npartitions\}}\stepsize}{\Npartitions} \circ \flow{\alpha_1^{\{\Npartitions-1\}}\stepsize}{\Npartitions-1} \circ \ldots \circ \flow{\alpha_1^{\{1\}}\stepsize}{1}
\end{equation}
for $\Npartitions \geq 2$ partitions and fractions $\alpha_i^{\{m\}}$ of the time step. 
The main idea is to iteratively employ splitting methods for two-split systems, which have been extensively studied in the literature, see e.g.~\cite{McLachlan_Quispel_2002,blanes2024splitting} and the references therein.

\subsection{Iterative splitting of the vector field}
Let $\rhs{\node{v}}(\x) \coloneqq \sum_{m \in \node{v}} \rhs{m}(\x)$ for any subset $v \subset \{1,\ldots,\Npartitions\}$.
We introduce the Lie derivatives $\LieDerivative{\node{v}} = \langle \rhs{\node{v}}(\cdot), \nabla_{\x} \rangle$ and $\LieDerivative{} = \langle \rhs{}(\cdot), \nabla_{\x} \rangle$, where $\langle \cdot, \cdot \rangle$ denotes the Euclidean inner product, i.e., $\langle \rhs{\node{v}}(\cdot), \nabla_{\x} \rangle = \sum_{k=1}^{n} \rhs{\node{v}}_k(\cdot)\, \partial_{x_k}$.
These differential operators allow us to formally express the $t$-flows of $\x^{\prime} = \rhs{\node{v}}(\x)$ and $\x^{\prime} = \rhs{}(\x)$ as $\flow{t}{\node{v}} = \exp(t \LieDerivative{\node{v}}) \mathrm{Id}$ and $\flow{t}{} = \exp(t \LieDerivative{}) \mathrm{Id}$, where $\mathrm{Id}$ denotes the identity map.
The iterative process starts by splitting the vector field $\rhs{} = \rhs{\node{r}}$ with $\node{r} \coloneqq \{1,\ldots,\Npartitions\}$ into two partitions
\begin{equation*}
    \begin{aligned}
        \rhs{\node{r}}(\x) &= \rhs{\node{i}}(\x) + \rhs{\node{j}}(\x), & \node{i},\node{j} &\subset \{1,\ldots,\Npartitions\}, & \node{i} \mathbin{\dot{\cup}} \node{j} &= \node{r},
    \end{aligned}
\end{equation*}
where $\mathbin{\dot{\cup}}$ denotes the disjoint union.
We compute a numerical approximation to the $t$-flow $\flow{t}{\node{r}} \coloneqq \flow{t}{}$ of \eqref{eq:split-system} with a splitting method for two-split systems,
$$\SplittingMethod{\stepsize}{\node{r}} = \flow{b_{s^{\{\node{r}\}}}^{\{\node{r}\}} \stepsize}{\node{j}} \circ \flow{a_{s^{\{\node{r}\}}}^{\{\node{r}\}} \stepsize}{\node{i}} \circ \ldots \circ \flow{b_1^{\{\node{r}\}} \stepsize}{\node{j}} \circ \flow{a_1^{\{\node{r}\}} \stepsize}{\node{i}}.$$
For an iterative splitting, we need to further partition $\rhs{\node{i}}$ and $\rhs{\node{j}}$. Let us consider the set $\node{i}$.
If $\lvert \node{i} \rvert = 1$, the elementary subsystem is not split further:  we either compute $\flow{}{\node{i}}$ exactly or approximate the flow with an appropriate numerical integration scheme $\SplittingMethod{}{\node{i}}$. 
If $\lvert \node{i} \rvert > 1$, we split the set into $\node{i} = \node{i}_{\mathrm{L}} \mathbin{\dot{\cup}} \node{i}_{\mathrm{R}}$ and apply another splitting method for two-split systems, which replaces $\flow{a_j^{\{\node{r}\}}\stepsize}{\node{i}}$ with
$$ \SplittingMethod{a_j^{\{\node{r}\}}\stepsize}{\node{i}} = \flow{a_j^{\{\node{r}\}} b_{s^{\{\node{i}\}}}^{\{\node{i}\}} \stepsize}{\node{i}_{\mathrm{R}}} \circ \flow{a_j^{\{\node{r}\}} a_{s^{\{\node{i}\}}}^{\{\node{i}\}} \stepsize}{\node{i}_{\mathrm{L}}} \circ \ldots \circ \flow{a_j^{\{\node{r}\}} b_{1}^{\{\node{i}\}} \stepsize}{\node{i}_{\mathrm{R}}} \circ \flow{a_j^{\{\node{r}\}} a_{1}^{\{\node{i}\}} \stepsize}{\node{i}_{\mathrm{L}}},\quad j=1,\ldots,s^{\{\node{r}\}}. $$
This process is iteratively repeated until a partition in singletons is achieved, and suitable numerical integration schemes are selected to compute the flows of the elementary subsystems.
This \emph{hierarchical splitting} process can be represented in terms of a \emph{splitting tree}.
\definition[Splitting tree]{
    A splitting tree for $\Npartitions \in \mathbb{N}$ partitions is a full ordered binary tree\footnote{%
    A full binary tree is a binary tree in which each node has either zero or two child nodes.
    An ordered binary tree is a rooted tree in which an ordering is specified for the children of each node. 
    In particular, we will distinguish between the left child node $\node{v}_{\mathrm{L}}$ and the right child node $\node{v}_{\mathrm{R}}$ of a node $\node{v} \in \nodeset$.
    } 
    $\graph=(\nodeset,\edges)$ with node set $\nodeset \subset \mathbb{P}_{\geq 1}(\{1,\ldots,\Npartitions\})$ so that the root node is $\node{r} =\{1,\ldots,\Npartitions\}$; if $\node{i}\!\in\! \nodeset$ has child nodes, then $\node{i}_{\mathrm{L}} \mathbin{\dot{\cup}} \node{i}_{\mathrm{R}}  =  \node{i}$; and each leaf node $\node{\ell} \!\in\! \nodeset$ has cardinality one, $\lvert \node{\ell}\rvert = 1$.
}\normalfont\medskip

\begin{figure}[htb]
    \centering
    \includegraphics{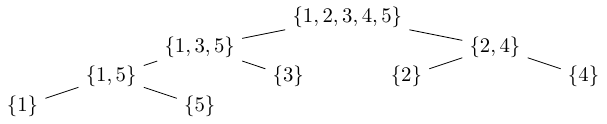}
    \caption{A splitting tree for $\Npartitions = 5$ partitions.}
    \label{fig:N5_splitting-tree}
\end{figure}

Consequently, a splitting tree has exactly $\Npartitions$ leaf nodes and satisfies $\lvert \nodeset\rvert = 2\Npartitions-1$.
An illustrative example of a splitting tree for $\Npartitions=5$ partitions is depicted in \Cref{fig:N5_splitting-tree}. 

Consider a splitting tree $\graph=(\nodeset,\edges)$. 
We define the sets $\innernodes \subset \nodeset$ of all inner nodes and $\leafnodes \subset \nodeset$ of all leaf nodes, $\nodeset = \innernodes \mathbin{\dot{\cup}} \leafnodes$. 
For any $\node{v} \in \nodeset$, let $\mathcal{P}^{\{\node{v}\}} \subset \innernodes$ denote the set of all inner nodes in the unique path from the root node $\node{r}$ to the node $\node{v}$, excluding $\node{v}$ itself. If $\node{v}$ is included, we write $\bar{\mathcal{P}}^{\{\node{v}\}} \coloneqq \mathcal{P}^{\{\node{v}\}} \cup \{\node{v}\}$.
\vspace*{-3ex}
\definition[Hierarchical splitting method]{
A hierarchical splitting method is a splitting method $\SplittingMethod{\stepsize}{}$ for $\Npartitions$-split systems \eqref{eq:split-system} that is defined by a splitting tree $\graph=(\nodeset,\edges)$ for $\Npartitions$ partitions such that (i) to each inner node $\node{i} \in \innernodes$, a splitting method for two-split systems is assigned, 
\begin{equation}\label{eq:HSM_inner-node_method}
\SplittingMethod{\stepsize}{\node{i}} = \flow{b_{s^{\{\node{i}\}}}^{\{\node{i}\}} \stepsize}{\node{i}_{\mathrm{R}}} \circ \flow{a_{s^{\{\node{i}\}}}^{\{\node{i}\}} \stepsize}{\node{i}_{\mathrm{L}}} \circ \ldots \circ \flow{b_1^{\{\node{i}\}} \stepsize}{\node{i}_{\mathrm{R}}} \circ \flow{a_1^{\{\node{i}\}} \stepsize}{\node{i}_{\mathrm{L}}} ;
\end{equation}
(ii) to each leaf node $\node{\ell} \in \leafnodes^{\prime} \subset \leafnodes$, a numerical integration scheme $\SplittingMethod{\stepsize}{\node{\ell}}$ consistent with $\x^{\prime} = \rhs{\node{\ell}}(\x)$ is assigned; (iii) to each leaf node $\node{\ell} \in \leafnodes \setminus \leafnodes^{\prime}$, its exact flow $\flow{h}{\node{\ell}}$ is assigned.
}\normalfont\medskip

The computational process of a hierarchical splitting method is illustrated in Algorithm~\ref{alg:HSM}.\medskip 

\begin{algorithm2e}[H]
\scriptsize
\caption{Hierarchical splitting method}
\Function{\textsc{Step}($\graph,\x,\stepsize,\node{v}$)\sidenote{$/\ast$ Applies the method assigned to node $\node{v}$ of the splitting tree $\graph$ to advance the state $\x$ by a step of size $\stepsize$.
    A full time step of the hierarchical splitting method is performed for $\node{v} = \node{r}$. $\ast/$}}{
    \If{$\node{v} \in \innernodes$}{

        \For{$j=1,\ldots,s^{\{\node{v}\}}$}{
            $\x \gets \textsc{Step}(\graph,\x,a_j^{\{\node{v}\}}\stepsize,\node{v}_{\mathrm{L}})$\;
            $\x \gets \textsc{Step}(\graph,\x,b_j^{\{\node{v}\}}\stepsize,\node{v}_{\mathrm{R}})$\;
        }

    }
    \ElseIf{$\node{v} \in \leafnodes^{\prime}$}{
        $\x \gets \SplittingMethod{\stepsize}{\node{v}}(\x)$\;
    }
    \Else{
        $\x \gets \flow{\stepsize}{\node{v}}(\x)$\;
    }

    \KwRet{$\x$}\;
}

\label{alg:HSM}
\end{algorithm2e}

\example
{\label{ex:hierarchical_splitting-method}
Consider the splitting tree from \Cref{fig:N5_splitting-tree}. Assign the Strang splitting \eqref{eq:Strang_N2} to all inner nodes, and the exact flow $\flow{\stepsize}{m}\ (m=1,\ldots,5)$ to all leaf nodes. 
The resulting hierarchical splitting method reads 
\begin{multline}\label{eq:example_HSM}
    \SplittingMethod{\stepsize}{} =  
    \flow{\stepsize/8}{1} \circ
    \flow{\stepsize/4}{5} \circ
    \flow{\stepsize/8}{1} \circ
    \flow{\stepsize/2}{3} \circ
    \flow{\stepsize/8}{1} \circ
    \flow{\stepsize/4}{5} \circ
    \flow{\stepsize/8}{1} \circ
    \flow{\stepsize/2}{2} \circ
    \flow{\stepsize}{4} \circ
    \flow{\stepsize/2}{2} \\
    \circ \flow{\stepsize/8}{1} \circ \flow{\stepsize/4}{5} \circ 
    \flow{\stepsize/8}{1} \circ
    \flow{\stepsize/2}{3} \circ
    \flow{\stepsize/8}{1} \circ
    \flow{\stepsize/4}{5} \circ
    \flow{\stepsize/8}{1}.
\end{multline}
}\normalfont

\subsection{Convergence order and accuracy}
The order conditions of splitting methods for $\Npartitions$-split systems can be derived from the Baker--Campbell--Hausdorff formula, or strictly speaking via its generalization to $\Npartitions \geq 2$ partitions \cite{SPITERI2026115271}. 
For investigating the convergence order and the leading order error terms of hierarchical splitting methods, we will make use of the following lemma.
\lemma{\label{lemma:BCH}
    Consider (non-commuting) matrices $\boldsymbol{A}^{\{m\}}, \boldsymbol{B}^{\{m\}} \in \mathbb{R}^{n \times n}$ and arbitrary scalar parameters $\alpha_i^{\{m\}},\beta_i^{\{m\}} \in \mathbb{R}$ for $i=1,\ldots,s$ with $s \in \mathbb{N}$, and $m=1,\ldots,\Npartitions$ with $\Npartitions \in \mathbb{N}$. 
    Additionally, let $\convergenceorder \in \mathbb{N}$. 
    Then, as $\stepsize \to 0$,  
    \begin{multline*}
        \prod_{i=1}^s \prod_{m=1}^{\Npartitions} \exp\big(\alpha_i^{\{m\}} \stepsize \boldsymbol{A}^{\{m\}} + \beta_i^{\{m\}} \stepsize^{\convergenceorder+1} \boldsymbol{B}^{\{m\}} + \mathcal{O}(\stepsize^{\convergenceorder+2})\big)
        \\
        = \exp\Big( \boldsymbol{C}(h) + \stepsize^{\convergenceorder+1} \Big(\sum_{m=1}^{\Npartitions} \sum_{i=1}^s  \beta_i^{\{m\}} \boldsymbol{B}^{\{m\}} \Big) + \mathcal{O}(\stepsize^{\convergenceorder+2})\Big),
    \end{multline*}
    where $\exp(\boldsymbol{C}(h)) = \left(\prod_{i=1}^s \prod_{m=1}^{\Npartitions} \exp\big(\alpha_i^{\{m\}} \stepsize \boldsymbol{A}^{\{m\}}\big)\right)$.
}
\begin{proof}
    A first-order Taylor expansion of the individual factors yields 
    \begin{multline}\label{eq:lemma1_eq}
        \exp\big(\alpha_i^{\{m\}} \stepsize \boldsymbol{A}^{\{m\}} + \beta_i^{\{m\}} \stepsize^{\convergenceorder+1} \boldsymbol{B}^{\{m\}} + \mathcal{O}(\stepsize^{\convergenceorder+2})\big) \\
        = \; \exp\big(\alpha_i^{\{m\}} \stepsize \boldsymbol{A}^{\{m\}}\big) + \beta_i^{\{m\}} \stepsize^{\convergenceorder+1} \boldsymbol{B}^{\{m\}} + \mathcal{O}(\stepsize^{\convergenceorder+2}).
    \end{multline}
    We can substitute \eqref{eq:lemma1_eq} into the product over $i$ and $m$ and expand. Since every perturbation $\beta_i^{\{m\}} \stepsize^{\convergenceorder+1}\boldsymbol{B}^{\{m\}} + \mathcal{O}(\stepsize^{\convergenceorder+2})$ is of order $\stepsize^{\convergenceorder + 1}$, any term involving a product of perturbations contributes to the remainder term $\mathcal{O}(\stepsize^{\convergenceorder + 2})$.
    Therefore, only the terms containing exactly one perturbation enter the leading order error term.
    Consequently, collecting all terms and inserting $\exp(\boldsymbol{C}(h))$ yields 
    $$ \exp(\boldsymbol{C}(h)) + \stepsize^{\convergenceorder+1} \Big(\sum_{m=1}^{\Npartitions} \sum_{i=1}^s  \beta_i^{\{m\}} \boldsymbol{B}^{\{m\}} \Big) + \mathcal{O}(\stepsize^{\convergenceorder+2}), $$
    which can be interpreted as a first-order Taylor expansion of the asserted expression.
\end{proof}

\theorem{\label{theorem:convergence_order}
Consider a hierarchical splitting method $\SplittingMethod{\stepsize}{}$ with numerical integration schemes $\SplittingMethod{\stepsize}{\node{v}}$ of order $\convergenceorder^{\{\node{v}\}}$ assigned to nodes $\node{v} \in \nodeset$, 
$$
\SplittingMethod{\stepsize}{\node{v}} = \exp\left( \stepsize \LieDerivative{\node{v}} + \stepsize^{\convergenceorder^{\{\node{v}\}}+1} \mathcal{O}_{\convergenceorder^{\{\node{v}\}}+1}^{\{\node{v}\}} + \mathcal{O}(\stepsize^{\convergenceorder^{\{\node{v}\}}+2})\right) \mathrm{Id},
$$
with leading error term $\mathcal{O}_{\convergenceorder^{\{\node{v}\}}+1}^{\{\node{v}\}}$.
For each node $\node{v} \in \nodeset$, we define the mapping
\begin{equation}\label{eq:Path01}
\begin{aligned}
    &\delta^{\{\node{v}\}} \colon \mathcal{P}^{\{\node{v}\}} \to \{0,1\}, &
    \delta^{\{\node{v}\}}(\node{w}) &= \begin{cases}
    1 , & \text{$\node{w}_{\mathrm{L}} \in \bar{\mathcal{P}}^{\{\node{v}\}}$,} \\
    0 , & \text{$\node{w}_{\mathrm{R}} \in \bar{\mathcal{P}}^{\{\node{v}\}}$.}
\end{cases} 
\end{aligned}
\end{equation} 
Then, $\SplittingMethod{\stepsize}{}$ is convergent of order at least $\convergenceorder \coloneqq \min_{\node{v} \in \innernodes \mathbin{\dot{\cup}} \leafnodes^{\prime}} \convergenceorder^{\{\node{v}\}}$, 
$$ \SplittingMethod{\stepsize}{} = \exp\left( \stepsize\LieDerivative{} + \stepsize^{\convergenceorder + 1}  \mathcal{O}_{\convergenceorder + 1} + \mathcal{O}(\stepsize^{\convergenceorder + 2}) \right)\mathrm{Id},$$
with leading error term
\begin{align}\label{eq:leading_error_term}
     \mathcal{O}_{\convergenceorder + 1} = &\sum_{\substack{\node{v} \in \innernodes \mathbin{\dot{\cup}} \leafnodes^{\prime} \\ \convergenceorder^{\{\node{v}\}} = \convergenceorder}} \left[ \prod_{\node{w} \in \mathcal{P}^{\{\node{v}\}}} \sum_{j=1}^{s^{\{\node{w}\}}} \Big( \delta^{\{\node{v}\}}(\node{w}) a_j^{\{\node{w}\}} + (1 - \delta^{\{\node{v}\}}(\node{w})) b_j^{\{\node{w}\}}\Big)^{\convergenceorder+1} \right] \mathcal{O}_{\convergenceorder + 1}^{\{\node{v}\}}.
\end{align} 
}
\begin{proof}
    For each node $\node{v} \in \nodeset$, the assigned method $\SplittingMethod{\stepsize}{\node{v}}$ or flow $\flow{\stepsize}{\node{v}}$ is employed $\prod_{\node{w} \in \mathcal{P}^{\{\node{v}\}}} s^{\{\node{w}\}}$ times with step sizes
    $$\tilde{\stepsize}^{\{\node{v}\}}_{(j_{\node{w}})_{\node{w} \in \mathcal{P}^{\{\node{v}\}}}} \coloneqq \prod_{\node{w} \in \mathcal{P}^{\{\node{v}\}}} \big(\delta^{\{\node{v}\}}(\node{w}) a_{j_{\node{w}}}^{\{\node{w}\}} + (1 - \delta^{\{\node{v}\}}(\node{w})) b_{j_{\node{w}}}^{\{\node{w}\}}\big)\stepsize,$$
    for $j_{\node{w}} = 1,\ldots,s^{\{\node{w}\}}$, and all $\node{w} \in \mathcal{P}^{\{\node{v}\}}$. 
    Let $\bm{H}^{\{\node{v}\}}$ denote the set of all step sizes $\tilde{\stepsize}^{\{\node{v}\}}_{(j_{\node{w}})_{\node{w} \in \mathcal{P}^{\{\node{v}\}}}}$.
    Note that $\sum_{\tilde{h} \in \bm{H}^{\{\node{v}\}}} \tilde{\stepsize} = \stepsize$ for all $\node{v} \in \nodeset$, provided that all splitting methods are consistent, i.e., $\sum_j a_j^{\{\node{i}\}} = \sum_j b_j^{\{\node{i}\}} = 1$ for all $\node{i} \in \innernodes$. 
    As per Gröbner's Lemma \cite[III Lemma 5.1]{HairerLubichWanner}, we can formally write $\SplittingMethod{\stepsize}{}$ as a product of 
    $ \exp(\tilde{\stepsize} \LieDerivative{\node{k}}), \, \node{k} \in \leafnodes \setminus \leafnodes^{\prime}, \, \tilde{h} \in \bm{H}^{\{\node{k}\}}$, and 
    $$ \exp\left( \tilde{\stepsize} \LieDerivative{\node{\ell}} + \tilde{\stepsize}^{\convergenceorder^{\{\node{\ell}\}} + 1} \mathcal{O}_{\convergenceorder^{\{\node{\ell}\}} + 1}^{\{\node{\ell}\}}  + \mathcal{O}\big(\tilde{\stepsize}^{\convergenceorder^{\{\node{\ell}\}} + 2}\big)\right), $$
    for $\node{\ell} \in \leafnodes^{\prime}$ and $\tilde{h} \in \bm{H}^{\{\node{\ell}\}} $.
    Applying Lemma \ref{lemma:BCH} for $\convergenceorder\coloneqq\min_{\node{\ell} \in \leafnodes^{\prime}} \convergenceorder^{\{\node{\ell}\}}$ yields
    \begin{equation}\label{eq:Theorem1_intermediate_exp}
        \SplittingMethod{\stepsize}{} = \exp\big( \stepsize\tilde{\LieDerivative{}} + \stepsize^{\convergenceorder + 1} \tilde{\mathcal{O}}_{\convergenceorder + 1} + \mathcal{O}(\stepsize^{\convergenceorder + 2}) \big) \mathrm{Id},
    \end{equation}
    where 
    \begin{equation*}
        \tilde{\mathcal{O}}_{\convergenceorder + 1} = \sum_{\substack{\node{\ell} \in \leafnodes^{\prime} \\ \convergenceorder^{\{\node{\ell}\}} = \convergenceorder}} \Bigg[ \prod_{\node{w} \in \mathcal{P}^{\{\node{\ell}\}}} \sum_{j=1}^{s^{\{\node{w}\}}} \Big( \delta^{\{\node{\ell}\}}(\node{w}) a_j^{\{\node{w}\}} + (1 - \delta^{\{\node{\ell}\}}(\node{w})) b_j^{\{\node{w}\}}\Big)^{\convergenceorder+1} \Bigg] \mathcal{O}_{\convergenceorder + 1}^{\{\node{\ell}\}},
    \end{equation*}
    and $\tilde{\SplittingMethod{}{}}_\stepsize = \exp(\stepsize \tilde{\LieDerivative{}}) \mathrm{Id}$ denotes the splitting method $\SplittingMethod{\stepsize}{}$ where all elementary subsystems are solved exactly, i.e., to each leaf node $\node{\ell} \in \leafnodes$ its exact flow $\flow{\stepsize}{\node{\ell}}$ is assigned.
    We now define the set $\boldsymbol{A} \subset \innernodes$ of previously investigated inner nodes, initially set to $\boldsymbol{A} = \emptyset$, and the set $\boldsymbol{B} \subset \nodeset \setminus \boldsymbol{A}$ of not yet investigated nodes whose child nodes are in $\boldsymbol{A} \mathbin{\dot{\cup}} \leafnodes$.
    Then, we can now interpret $\tilde{\SplittingMethod{}{}}_{\stepsize} = \exp(\stepsize\tilde{\LieDerivative{}})\mathrm{Id}$ as the splitting method $\SplittingMethod{\stepsize}{}$ where all methods assigned to nodes $\node{v} \in \boldsymbol{A} \mathbin{\dot{\cup}}\leafnodes^{\prime}$ are replaced with the exact flow $\flow{\stepsize}{\node{v}} = \exp(\stepsize \LieDerivative{\node{v}})\mathrm{Id}$ and write \eqref{eq:Theorem1_intermediate_exp} with
    \begin{equation}\label{eq:Theorem1-intermediate} 
    \tilde{\mathcal{O}}_{\convergenceorder + 1} = \sum_{\substack{\node{v} \in \boldsymbol{A} \mathbin{\dot{\cup}} \leafnodes^{\prime} \\ \convergenceorder^{\{\node{v}\}} = \convergenceorder}} \Bigg[ \prod_{\node{w} \in \mathcal{P}^{\{\node{v}\}}} \sum_{j=1}^{s^{\{\node{w}\}}} \Big( \delta^{\{\node{v}\}}(\node{w}) a_j^{\{\node{w}\}} + (1 - \delta^{\{\node{v}\}}(\node{w})) b_j^{\{\node{w}\}}\Big)^{\convergenceorder+1} \Bigg] \mathcal{O}_{\convergenceorder + 1}^{\{\node{v}\}} , 
    \end{equation}
    and $\convergenceorder \coloneqq \min_{\node{v} \in \boldsymbol{A} \mathbin{\dot{\cup}} \leafnodes^{\prime}} \convergenceorder^{\{\node{v}\}}$.
    As long as $\boldsymbol{A} \neq \innernodes$, we proceed as follows.
    {
    \setlength{\leftmargini}{1.5em}
    \begin{enumerate}
        \item[\bfseries 1)] Select a node $\node{i} \in \boldsymbol{B}$ and consider $\node{i}_{\mathrm{L}}, \node{i}_{\mathrm{R}} \in \boldsymbol{A} \mathbin{\dot{\cup}} \leafnodes$.
        \item[\bfseries 2)] The method $\tilde{\SplittingMethod{}{}}_\stepsize$ involves applications of the splitting method 
        $$ \SplittingMethod{\tilde{\stepsize}}{\node{i}} = \exp\!\big( a_1^{\{\node{i}\}} \tilde{\stepsize} \LieDerivative{\node{i}_{\mathrm{L}}} \big) \exp\!\big( b_1^{\{\node{i}\}} \tilde{\stepsize} \LieDerivative{\node{i}_{\mathrm{R}}} \big) \cdots \exp\!\big( a_{s^{\{\node{i}\}}}^{\{\node{i}\}} \tilde{\stepsize} \LieDerivative{\node{i}_{\mathrm{L}}} \big) \exp\!\big( b_{s^{\{\node{i}\}}}^{\{\node{i}\}} \tilde{\stepsize} \LieDerivative{\node{i}_{\mathrm{R}}} \big) \mathrm{Id}$$
        with $\tilde{\stepsize} \in \bm{H}^{\{\node{i}\}}$. As per assumption, we can replace these occurrences with
        $$ \SplittingMethod{\tilde{\stepsize}}{\node{i}} = \exp\big( \tilde{\stepsize} \LieDerivative{\node{i}} + \tilde{\stepsize}^{\convergenceorder^{\{\node{i}\}} + 1} \mathcal{O}_{\convergenceorder^{\{\node{i}\}} + 1}^{\{\node{i}\}} + \mathcal{O}(\tilde{\stepsize}^{\convergenceorder^{\{\node{i}\}} + 2}) \big) \mathrm{Id}. $$
        \item[\bfseries 3)] Add $\node{i}$ to the set $\boldsymbol{A}$ and, if $\node{i} \neq \node{r}$, add the parent node of $\node{i}$ to the set $\boldsymbol{B}$. 
        Subsequently, formally applying Lemma~\ref{lemma:BCH} with $\convergenceorder \coloneqq \min_{\node{v} \in \boldsymbol{A} \mathbin{\dot{\cup}} \leafnodes^{\prime}} \convergenceorder^{\{\node{v}\}}$ once more yields \eqref{eq:Theorem1-intermediate}.
    \end{enumerate}
    }
    Repeating this process until $\boldsymbol{A} = \innernodes$, we have $\boldsymbol{A} \mathbin{\dot{\cup}} \leafnodes^{\prime} = \innernodes \mathbin{\dot{\cup}} \leafnodes^{\prime}$, $\tilde{\SplittingMethod{}{}}_\stepsize = \flow{\stepsize}{} = \exp(\stepsize \LieDerivative{})\mathrm{Id}$, and thus obtain the asserted expression.
\end{proof}\normalfont
\corollary{
As per \Cref{theorem:convergence_order}, assigning splitting methods and numerical integration schemes of order at least $\convergenceorder$ to all nodes within the splitting tree is a sufficient condition to obtain a hierarchical splitting method of order $\convergenceorder$.}\medskip\normalfont 

This straightforward extension of the order theory facilitates the derivation of hierarchical splitting methods of order $\convergenceorder>2$ by utilizing efficient splitting methods for two-split systems that are already available in the literature, such as the methods mentioned in the survey papers \cite{McLachlan_Quispel_2002,blanes2024splitting} on splitting methods and the references therein.

\remark
{ 
Methods for $(K>2)$-split systems are also available in the literature~\cite{auzinger2015splitting,auzinger2017splitting,spiteri_beyond_2025,SPITERI2026115271}. 
By generalizing the concept of a splitting tree, these methods can be integrated into the hierarchical splitting approach:
instead of requiring a splitting tree to be a binary tree, we require that an inner node $\node{i} \in \innernodes$ to which a splitting method for $K$-split systems is assigned has exactly $K$ child nodes.
}\normalfont

\proposition[Self-adjoint hierarchical splitting methods]{
A hierarchical splitting method $\SplittingMethod{\stepsize}{}$ is self-adjoint, i.e., $\SplittingMethod{\stepsize}{} = \SplittingMethod{-\stepsize}{}^{-1} \eqqcolon \SplittingMethod{\stepsize}{}^*$, provided that the methods assigned to the leaf nodes are self-adjoint and the splitting methods~\eqref{eq:HSM_inner-node_method} assigned to the inner nodes have palindromic sequences of sub-steps.
}
\begin{proof}
    Recall that $(\SplittingMethod{\stepsize}{1} \circ \SplittingMethod{\stepsize}{2})^* = (\SplittingMethod{\stepsize}{2})^* \circ (\SplittingMethod{\stepsize}{1})^*$ and that self-adjointness is preserved under scaling of the step size.
    We proceed by structural induction over the splitting tree.
    For the leaf nodes $\node{\ell} \in \leafnodes$, the assigned (numerical) flows are self-adjoint as per assumption.
    Now let $\node{i} \in \innernodes$ and assume that the maps assigned to $\node{i}_{\mathrm{L}}$ and $\node{i}_{\mathrm{R}}$ are self-adjoint.
    Taking the adjoint of the splitting method \eqref{eq:HSM_inner-node_method} reverses the sequence of sub-steps and replaces each of them by its adjoint, which leaves them unchanged by the induction hypothesis.
    By assumption, the sequence of sub-steps of the splitting method assigned to node $\node{i}$ is palindromic, so that the reversed sequence coincides with the original one and $\SplittingMethod{\stepsize}{\node{i}}$ is self-adjoint.
    For $\node{i} = \node{r}$, this yields the self-adjointness of the hierarchical splitting method.
\end{proof}
\normalfont

\remark{\label[remark]{remark:commutators}
Besides conventional splitting methods for two-split systems, it is also feasible to assign commutator-based splitting methods (cf.~\cite{wisdom1996symplectic,LopezMarcos1997HessianVectorProducts,omelyan2003symplectic,schäfers2024hessianfree,MONCH202525}) to nodes of a splitting tree. 
Assigning a commutator-based splitting method to an inner node $\node{i} \in \innernodes$ requires that the approach is consistent with the Lie derivatives of the child nodes $\node{i}_{\mathrm{L}}$ and $\node{i}_{\mathrm{R}}$.
}

\section{Multirate hierarchical splitting methods}
\label{sec:multiple_time_stepping}
When considering $\Npartitions$-split systems, it is frequently the case that computationally expensive parts of the right-hand side do not contribute much to the dynamics of the solution, suggesting the use of \emph{multiple time stepping techniques} that evaluate the expensive part of the vector field less often.
Multiple time stepping methods can be constructed within the framework of splitting methods \cite{BIESIADECKI1993}.
A common example is a multirate generalization of the Strang splitting that is known as the \emph{impulse method} \cite{HairerLubichWanner}.
It is given by
\begin{equation}\label{eq:impulse_method}
    \SplittingMethod{\stepsize}{} = \SplittingMethod{\stepsize/2}{\mathfrak{s}} \circ \big( \SplittingMethod{\stepsize/\mrfactor{}}{\mathfrak{f}} \big)^{\mrfactor{}} \circ \SplittingMethod{\stepsize/2}{\mathfrak{s}}, \quad \mrfactor{} \in \mathbb{N},
\end{equation}
where $\SplittingMethod{\stepsize}{\mathfrak{f}}$ and $\SplittingMethod{\stepsize}{\mathfrak{s}}$ are numerical integration schemes consistent with the $\mathfrak{f}$ast subsystem $\x^{\prime} = \rhs{\mathfrak{f}}(\x)$ and $\mathfrak{s}$low subsystem $\x^{\prime} = \rhs{\mathfrak{s}}(\x)$, respectively.
It was already noted in \cite{HairerLubichWanner} that multiple time stepping techniques can be iteratively extended to problems with $\Npartitions > 2$ partitions with different time scales. 
The idea is to further split the fast partition into $\rhs{\mathfrak{f}}(\x) = \rhs{\mathfrak{f,f}}(\x) + \rhs{\mathfrak{f,s}}(\x)$, and to replace the numerical integration scheme $\SplittingMethod{\stepsize}{\mathfrak{f}}$ with another multiple time stepping method.
This approach is widely used in molecular dynamics simulations~\cite{tuckerman_reversible_1992,SEXTON1992665,URBACH200687} and aligns with the splitting process of a hierarchical splitting method, which we formalize by introducing multirate splitting trees.
\definition[Multirate splitting tree]{
A multirate splitting tree is a splitting tree $\graph = (\nodeset,\edges)$ with weighted edges where for any node $\node{v} \in \nodeset \setminus \{\node{r}\}$, the respective multirate factor $\mrfactor{\node{v}} \in \mathbb{N}$ is assigned to the edge connecting the node $\node{v}$ with its parent node.
}\medskip\normalfont

Based on a multirate splitting tree, a multirate version is obtained by extending \eqref{eq:HSM_inner-node_method} to
\begin{multline}\label{eq:mr_naive_approach}
    \SplittingMethod{\stepsize}{\node{i}} = \Big( \flow{  b_{s^{\{\node{i}\}}}^{\{\node{i}\}} \stepsize\big/\mrfactor{\node{i}_{\mathrm{R}}} }{\node{i}_{\mathrm{R}}} \Big)^{\mrfactor{\node{i}_{\mathrm{R}}}} \circ \Big(\flow{a_{s^{\{\node{i}\}}}^{\{\node{i}\}} \stepsize\big/\mrfactor{\node{i}_{\mathrm{L}}} }{\node{i}_{\mathrm{L}}}\Big)^{\mrfactor{\node{i}_{\mathrm{L}}}} 
    \circ \ldots \\
    \circ \Big(\flow{b_{1}^{\{\node{i}\}} \stepsize\big/\mrfactor{\node{i}_{\mathrm{R}}} }{\node{i}_{\mathrm{R}}}\Big)^{\mrfactor{\node{i}_{\mathrm{R}}}} \circ \Big(\flow{a_1^{\{\node{i}\}} \stepsize\big/\mrfactor{\node{i}_{\mathrm{L}}} }{\node{i}_{\mathrm{L}}}\Big)^{\mrfactor{\node{i}_{\mathrm{L}}}}.
\end{multline}
In fact, \eqref{eq:mr_naive_approach} will introduce a modification only if numerical flow approximations are used. For any analytically computed flow assigned to a leaf node $\node{\ell} \in \leafnodes \setminus \leafnodes^\prime$, we can select $\mrfactor{\node{\ell}} = 1$.

\begin{figure}[htb!]
    \begin{center}
    \includegraphics{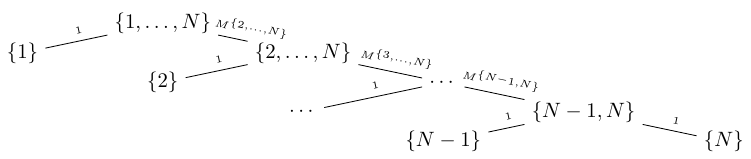}
\end{center}
    \caption{Multirate splitting tree of a generalized impulse method for $\Npartitions$-split systems.}
    \label{fig:mr_HSM_splitting_tree}
\end{figure}

\example{
The generalization of the impulse method \eqref{eq:impulse_method} to $\Npartitions > 2$ partitions can be represented by a multirate hierarchical splitting method where the Strang splitting \eqref{eq:Strang_N2} is assigned to all inner nodes $\node{i} \in \innernodes$ of the multirate splitting tree in \Cref{fig:mr_HSM_splitting_tree}.
}

Typically, we determine the multirate factor $\mrfactor{} \in \mathbb{N}$ such that the step size $\stepsize / \mrfactor{}$ is sufficiently small to meet the stability and accuracy requirements. 
For any node $\node{v} \in \nodeset \setminus \{\node{r}\}$, the assigned method is employed with step sizes 
\begin{equation}\label{eq:mr-approach_step_sizes}
    \frac{\prod_{\node{w} \in \mathcal{P}^{\{\node{v}\}}} \delta^{\{\node{v}\}}(\node{w}) a_{j_{\node{w}}}^{\{\node{w}\}} + (1-\delta^{\{\node{v}\}}(\node{w})) b_{j_{\node{w}}}^{\{\node{w}\}}}{ \prod_{\node{u} \in \bar{\mathcal{P}}^{\{\node{v}\}} \setminus \{\node{r}\} } \mrfactor{\node{u}} },
\end{equation}
where $j_{\node{w}} = 1,\ldots,s^{\{\node{w}\}}$ for all $\node{w} \in \mathcal{P}^{\{\node{v}\}}$, and $\bar{\mathcal{P}}^{\{\node{v}\}} \coloneqq \mathcal{P}^{\{\node{v}\}} \cup \{\node{v}\}$.
These step sizes can be vastly different in magnitude. 
For instance, the coefficients of the fourth-order splitting method of RKN type \eqref{eq:OMF4} satisfy $\max\nolimits_{i} \lvert a_i\rvert / \min\nolimits_j \lvert a_j \rvert \approx 17.23$, and $\max\nolimits_{i} \lvert b_i\rvert / \min\nolimits_j \lvert b_j \rvert \approx 8.12$.
The requirements for the multirate factor $\mrfactor{\node{v}}$ are driven by the step size \eqref{eq:mr-approach_step_sizes} of the largest modulus to sufficiently resolve the dynamics of the subsystem $\x^{\prime} = \rhs{\node{v}}(\x)$. 
For the sub-steps with step sizes \eqref{eq:mr-approach_step_sizes} of smaller modulus, a smaller multirate factor $\mrfactor{\node{v}}$ is already sufficient to meet the accuracy and stability requirements.
As a result, the approach \eqref{eq:mr_naive_approach} may cause computational overhead. This overhead can be mitigated by reweighting the multirate factors for each individual sub-step.
Suppose a sub-step requires advancing the solution by $c \cdot \stepsize$, $c \in \mathbb{R} \setminus\{0\}$, and the desired step size is $\stepsize/\mrfactor{}$, $\mrfactor{} \in \mathbb{N}$. 
We aim to ensure that each sub-step uses a step size that is less than or equal to the desired step size $\stepsize/\mrfactor{}$. 
At the same time, the step size should remain as close to $\stepsize/\mrfactor{}$ as possible to minimize computational overhead.

\lemma{\label{lemma:mr-factor_optimality}
The smallest $\tilde{\mrfactor{}} \in \mathbb{N}$ satisfying $\lvert c \cdot \stepsize \rvert / \tilde{\mrfactor{}} \leq \stepsize / \mrfactor{}$ for a given constant $c \in \mathbb{R} \setminus \{0\}$, step size $\stepsize > 0$, and multirate factor $\mrfactor{} \in \mathbb{N}$ is
\begin{equation}\label{eq:reweighted_mr-factor}
    \tilde{\mrfactor{}} \coloneqq \lceil \lvert c \rvert \cdot \mrfactor{} \rceil.
\end{equation}}
\vspace{-2em}
\begin{proof}
    The choice \eqref{eq:reweighted_mr-factor} is feasible since 
    $$ \left\lvert \frac{c \cdot \stepsize}{\tilde{\mrfactor{}}} \right\rvert = \frac{\lvert c \rvert \cdot \stepsize}{\lceil \lvert c \rvert \cdot \mrfactor{} \rceil} \leq \frac{\lvert c \rvert \cdot \stepsize}{\lvert c \rvert \cdot \mrfactor{}} =\frac{\stepsize}{\mrfactor{}}. $$
    For $\tilde{\mrfactor{}} = 1$, optimality is evident. For $\tilde{\mrfactor{}} > 1$, we obtain the optimality of the choice \eqref{eq:reweighted_mr-factor} as
    \begin{equation*} 
        \frac{\stepsize}{\mrfactor{}} = \frac{\lvert c \rvert \cdot \stepsize}{\lvert c \rvert \cdot \mrfactor{}} \leq \frac{\lvert c \rvert \cdot \stepsize}{\lfloor \lvert c \rvert \cdot \mrfactor{} \rfloor} . \qedhere
    \end{equation*}
\end{proof}\normalfont

In particular, when computing a sub-step $\SplittingMethod{c\stepsize}{\node{v}}$ in node $\node{v}$, the coefficient $c$ reads 
\begin{equation}
c^{\{\node{v}\}}_{(j_{\node{w}})_{\node{w} \in \mathcal{P}^{\{\node{v}\}}}} = 
\frac{
    \prod_{\node{w} \in \mathcal{P}^{\{\node{v}\}}} \delta^{\{\node{v}\}}(\node{w}) a_{j_{\node{w}}}^{\{\node{w}\}} + (1-\delta^{\{\node{v}\}}(\node{w}))  b_{j_{\node{w}}}^{\{\node{w}\}} 
    }{ 
    \prod_{\node{u} \in \mathcal{P}^{\{\node{v}\}} \setminus \{\node{r}\}  } \mrfactor{\node{u}}_{(j_{\node{k}})_{\node{k} \in \mathcal{P}^{\{\node{u}\}} }} 
    }\, ,
\end{equation}
where $\mrfactor{\node{u}}_{(j_{\node{k}})_{\node{k} \in \mathcal{P}^{\{\node{u}\}}}}$ denotes the reweighted multirate factors of the nodes in $\mathcal{P}^{\{\node{v}\}}$. 
We define 
$$\mrfactor{\node{v}}_{(j_{\node{w}})_{\node{w} \in \mathcal{P}^{\{\node{v}\}}}} \coloneqq \Big\lceil \lvert c^{\{\node{v}\}}_{(j_{\node{w}})_{\node{w} \in \mathcal{P}^{\{\node{v}\}}}} \rvert \cdot \mrfactor{\node{v}} \Big\rceil.$$
For brevity, we write $c^{\{\node{v}\}}_{(j_{\node{w}})}$ and $\mrfactor{\node{v}}_{(j_{\node{w}})}$ instead of $c^{\{\node{v}\}}_{(j_{\node{w}})_{\node{w} \in \mathcal{P}^{\{\node{v}\}}}}$ and $\mrfactor{\node{v}}_{(j_{\node{w}})_{\node{w} \in \mathcal{P}^{\{\node{v}\}}}}$, respectively.

\definition[Multirate hierarchical splitting method]{
A multirate hierarchical splitting method is a splitting method $\SplittingMethod{\stepsize}{}$ for $\Npartitions$-split systems that is defined by a multirate splitting tree $\graph = (\nodeset,\edges)$ for $\Npartitions$ partitions such that 
{
\setlength{\leftmargini}{2em}
\begin{itemize}
    \item to each inner node $\node{i} \in \innernodes$, there is an assigned splitting method $\SplittingMethod{\stepsize}{\node{i}}$ for two-split systems. 
    In particular, if a sub-step advances the solution of the subsystem $\x^{\prime} = \rhs{\node{i}}(\x)$ by a step of size $c \cdot \stepsize$ with $c \in \mathbb{R} \setminus \{0\}$, then the sub-step reads 
    \begin{multline*}
        \SplittingMethod{c\stepsize}{\node{i}} = \Big( \flow{ \left. c b_{s^{\{\node{i}\}}}^{\{\node{i}\}} \stepsize\middle/\mrfactor{\node{i}_{\mathrm{R}}}_{s^{\{\node{i}\}}} \right. }{\node{i}_{\mathrm{R}}} \Big)^{\mrfactor{\node{i}_{\mathrm{R}}}_{s^{\{\node{i}\}}}} \circ \Big(\flow{\left.c a_{s^{\{\node{i}\}}}^{\{\node{i}\}} \stepsize\middle/\mrfactor{\node{i}_{\mathrm{L}}}_{s^{\{\node{i}\}}}\right. }{\node{i}_{\mathrm{L}}}\Big)^{\mrfactor{\node{i}_{\mathrm{L}}}_{s^{\{\node{i}\}}}} 
        \circ \ldots \\
        \circ \Big(\flow{\left. c b_{1}^{\{\node{i}\}} \stepsize\middle/\mrfactor{\node{i}_{\mathrm{R}}}_{1}\right. }{\node{i}_{\mathrm{R}}}\Big)^{\mrfactor{\node{i}_{\mathrm{R}}}_{1}} \circ \Big(\flow{\left.c a_1^{\{\node{i}\}} \stepsize\middle/\mrfactor{\node{i}_{\mathrm{L}}}_{1} \right.}{\node{i}_{\mathrm{L}}}\Big)^{\mrfactor{\node{i}_{\mathrm{L}}}_{1}}, 
    \end{multline*}
    with $\mrfactor{\node{i}_{\mathrm{L}}}_j \coloneqq \big\lceil \lvert c a_j^{\{\node{i}\}} \rvert \cdot \mrfactor{\node{i}_{\mathrm{L}}} \big\rceil $ and $\mrfactor{\node{i}_{\mathrm{R}}}_j \coloneqq \big\lceil \lvert c b_j^{\{\node{i}\}} \rvert \cdot \mrfactor{\node{i}_{\mathrm{R}}} \big\rceil$ for $ j=1,\ldots,s^{\{\node{i}\}}$;
    \item to each leaf node $\node{\ell} \in \leafnodes^{\prime}$, a numerical integration scheme $\SplittingMethod{\stepsize}{\node{\ell}}$ is assigned;
    \item to each leaf node $\node{\ell} \in \leafnodes \setminus \leafnodes^{\prime}$, its exact flow $\flow{\stepsize}{\node{\ell}}$ is assigned, the multirate factor is $\mrfactor{\node{\ell}} = 1$, and no reweighting is applied.
\end{itemize}
}
}\normalfont\medskip

The computational process of a multirate hierarchical splitting method is illustrated in Algorithm \ref{alg:MR_HSM}. 
The multiple time stepping approach evidently affects the truncation error of the numerical integration scheme. 
The qualitative statement of the convergence order, however, remains unchanged.

\begin{algorithm2e}[htb!]
\scriptsize
\caption{Multirate hierarchical splitting method}

\Function{\textsc{Step}($\graph,\x,c,\stepsize,\node{v}$)\sidenote{$/\ast$ Applies the method assigned to node $\node{v}$ of the multirate splitting tree $\graph$ to advance the state $\x$ by a step of size $c\stepsize$. A full time step of the multirate hierarchical splitting method is performed for $c=1$ and $\node{v} = \node{r}$. $\ast/$}}{
    \If{$\node{v} = \node{r}$}{
        $\tilde{\mrfactor{}} \gets 1$\;
    }
    \Else{
        $\tilde{\mrfactor{}} \gets \lceil \lvert c\rvert \cdot \mrfactor{\node{v}} \rceil$\; 
    }
    \If{$\node{v} \in \innernodes$}{
        \For{$m=1,\ldots,\tilde{\mrfactor{}}$}{
            \For{$j=1,\ldots,s^{\{\node{v}\}}$}{
                $\x \gets \textsc{Step}(\graph,\x,a_j^{\{\node{v}\}}c/\tilde{\mrfactor{}},\stepsize,\node{v}_{\mathrm{L}})$\;
                $\x \gets \textsc{Step}(\graph,\x,b_j^{\{\node{v}\}}c/\tilde{\mrfactor{}},\stepsize,\node{v}_{\mathrm{R}})$\;
            }
        }
    }
    \ElseIf{$\node{v} \in \leafnodes^{\prime}$}{
        \For{$m=1,\ldots,\tilde{\mrfactor{}}$}{
            $\x \gets \SplittingMethod{c\stepsize/\tilde{\mrfactor{}}}{\node{v}}(\x)$\;
        }
    }
    \Else{
        $\x \gets \flow{c\stepsize}{\node{v}}(\x)$\;
    }
    \KwRet{$\x$}\;
}
\label{alg:MR_HSM}
\end{algorithm2e}

\theorem{\label{theorem:MR_HSM}
Consider a multirate hierarchical splitting method $\SplittingMethod{\stepsize}{}$ with numerical integration schemes $\SplittingMethod{\stepsize}{\node{v}}$ of order $\convergenceorder^{\{\node{v}\}}$ assigned to nodes $\node{v} \in \nodeset$, 
$$\SplittingMethod{\stepsize}{\node{v}} = \exp\left( \stepsize \LieDerivative{\node{v}} + \stepsize^{\convergenceorder^{\{\node{v}\}} + 1} \mathcal{O}_{\convergenceorder^{\{\node{v}\}} + 1}^{\{\node{v}\}} + \mathcal{O}(\stepsize^{\convergenceorder^{\{\node{v}\}} + 2}) \right) \mathrm{Id},$$
with leading error terms $\mathcal{O}_{\convergenceorder^{\{\node{v}\}} + 1}^{\{\node{v}\}}$ and let $\delta^{\{\node{v}\}}$ denote the mappings defined in \eqref{eq:Path01}.
Then, $\SplittingMethod{\stepsize}{}$ is convergent of order at least $\convergenceorder \coloneqq \min_{\node{v} \in \innernodes \mathbin{\dot{\cup}} \leafnodes^{\prime}} \convergenceorder^{\{\node{v}\}}$, 
$$ \SplittingMethod{\stepsize}{} = \exp\left( \stepsize \LieDerivative{} + \stepsize^{\convergenceorder + 1}  \mathcal{O}_{\convergenceorder + 1} + \mathcal{O}(\stepsize^{\convergenceorder + 2}) \right),$$
with leading error term
\begin{align*}
     \mathcal{O}_{\convergenceorder + 1} 
     =
     \sum_{\substack{\node{v} \in \innernodes \mathbin{\dot{\cup}} \leafnodes^{\prime} \\ \convergenceorder^{\{\node{v}\}} = \convergenceorder}} 
     \left[  
     \sum_{(j_{\node{w}}) \in \mathcal{J}^{\{\node{v}\}}} 
     \frac{
        \left(\prod_{\node{w} \in \mathcal{P}^{\{\node{v}\}}} \delta^{\{\node{v}\}}(\node{w}) a_{j_{\node{w}}}^{\{\node{w}\}} + (1 - \delta^{\{\node{v}\}}(\node{w})) b_{j_{\node{w}}}^{\{\node{w}\}}\right)^{\convergenceorder + 1}
     }{
        \left(\prod_{\node{w} \in \bar{\mathcal{P}}^{\{\node{v}\}} \setminus \{\node{r}\} }
        \mrfactor{\node{w}}_{(j_{\node{u}})}
        \right)^\convergenceorder
     } 
     \right] 
     \mathcal{O}_{\convergenceorder + 1}^{\{\node{v}\}},
\end{align*} 
where $\mathcal{J}^{\{\node{v}\}} \coloneqq \left\{ (j_{\node{w}})_{\node{w} \in \mathcal{P}^{\{\node{v}\}}} \, : \, 1 \leq j_{\node{w}} \leq s^{\{\node{w}\}} \, \text{for all } \node{w} \in \mathcal{P}^{\{\node{v}\}} \right\}$.
}
\begin{proof}
    The root node $\node{r}$ is not impacted by the multiple time stepping approach as the assigned splitting method $\SplittingMethod{\stepsize}{\node{r}}$ still computes a single time step of step size $\stepsize$. 
    Additionally, all nodes $\node{\ell} \in \leafnodes \setminus \leafnodes^{\prime}$ remain unaffected, as the exact flow $\flow{\stepsize}{\node{\ell}}$ is assigned to these nodes.
    Consider any node $\node{v} \in \innernodes \mathbin{\dot{\cup}} \leafnodes^{\prime}$ where a numerical integration scheme $\SplittingMethod{\stepsize}{\node{v}}$ is assigned.
    At any occurrence, identified by a tuple $(j_{\node{w}}) \in \mathcal{J}^{\{\node{v}\}}$, the method computes $\mrfactor{\node{v}}_{(j_{\node{w}})} $ steps of step size $c_{(j_{\node{w}})}^{\{\node{v}\}} \stepsize \big/ \mrfactor{\node{v}}_{(j_{\node{w}})}$.  
    These consecutive applications can be formally expressed as
    \begin{align*}
        \bm{\Psi}_{c_{(j_{\node{w}})}^{\{\node{v}\}} \stepsize}^{\{\node{v}\}} &\coloneqq \left(\SplittingMethod{c_{(j_{\node{w}})}^{\{\node{v}\}} \stepsize\big/\mrfactor{\node{v}}_{(j_{\node{w}})}}{\node{v}}\right)^{\mrfactor{\node{v}}_{(j_{\node{w}})}} \\
        &= \exp \left( c_{(j_{\node{w}})}^{\{\node{v}\}} \stepsize \LieDerivative{\node{v}} + \mrfactor{\node{v}}_{(j_{\node{w}})}  \left(\frac{c_{(j_{\node{w}})}^{\{\node{v}\}} \stepsize} {\mrfactor{\node{v}}_{(j_{\node{w}})}}\right)^{\convergenceorder^{\{\node{v}\}}+1}
    \mathcal{O}_{\convergenceorder^{\{\node{v}\}} + 1}^{\{\node{v}\}} + \mathcal{O}\big(\stepsize^{\convergenceorder^{\{\node{v}\}} + 2}\big) \right) \mathrm{Id}.
    \end{align*}
    We can extend the iterative process outlined in the proof of \Cref{theorem:convergence_order} as follows: between steps 1) and 2), we replace the methods $\SplittingMethod{\stepsize}{\node{v}}$ with $\bm{\Psi}_{\stepsize}^{\{\node{v}\}}$ for all $\node{v} \in \nodeset \setminus \{\node{r}\}$.  
    In the proof of~\Cref{theorem:convergence_order}, we used in \eqref{eq:leading_error_term} that the sum and the product were interchangeable. 
    As the reweighted multirate factors are generally different, we can no longer interchange the sum and product signs.
    Consequently, the leading error term of the intermediate steps \eqref{eq:Theorem1-intermediate} becomes
    $$ \tilde{\mathcal{O}}_{\convergenceorder + 1} = \sum_{\substack{\node{v} \in \boldsymbol{A} \mathbin{\dot{\cup}} \leafnodes^{\prime} \\ \convergenceorder^{\{\node{v}\}} = \convergenceorder}}\left[ \sum_{(j_{\node{w}}) \in \mathcal{J}^{\{\node{v}\}}} 
    \frac{
        \left(\prod_{\node{w} \in \mathcal{P}^{\{\node{v}\}}} \delta^{\{\node{v}\}}(\node{w}) a_{j_{\node{w}}}^{\{\node{w}\}} + (1 - \delta^{\{\node{v}\}}(\node{w})) b_{j_{\node{w}}}^{\{\node{w}\}}\right)^{\convergenceorder + 1}
     }{
        \left(\prod_{\node{w} \in \bar{\mathcal{P}}^{\{\node{v}\}} \setminus \{\node{r}\} }
        \mrfactor{\node{w}}_{(j_{\node{u}})}
        \right)^\convergenceorder
     } 
    \right] \mathcal{O}_{\convergenceorder + 1}^{\{\node{v}\}} . $$
    The remainder of this proof is identical to the proof of \Cref{theorem:convergence_order}.
\end{proof}\normalfont

\remark{
As the number of applications of a method $\SplittingMethod{\stepsize}{\node{v}}$ grows exponentially with the distance of $\node{v}$ to $\node{r}$, a general guideline is to assign elementary subsystems $\x^\prime = \rhs{m}(\x)$ that require a smaller step size to leaf nodes with a greater distance to the root node. 
Consequently, a natural assumption is that the multirate factors along a path from the root node $\node{r}$ to a node $\node{v}$ are monotonically increasing. \medskip
}

In the limit $\mrfactor{\node{i}} \to \infty$, the subsystem $\x^\prime = \rhs{\node{i}}(\x)$ will be integrated exactly for any node $\node{i} \in \innernodes \setminus \{\node{r}\}$, provided that all methods assigned to nodes in this branch of the multirate splitting tree are consistent.
Therefore, we can assign lower-order schemes to nodes of the sub-tree without altering the convergence order of the overall method.
Since the limit $\mrfactor{\node{i}} \to \infty$ is not attained in practice, we impose conditions on the multirate factors ensuring that the scaling behavior for step sizes $h \geq h_{\min}$ is superior to the expected scaling based on the order theory.
We refer to this as the \emph{computational order} of a multirate integration scheme. 
The conditions on the multirate factors are derived in the following lemma.

\lemma{\label{lemma:computational_order}
Let $\stepsize > 0$, $c \in \mathbb{R} \setminus \{0\}$, $\mrfactor{} \in \mathbb{N}$, and $\convergenceorder,q \in \mathbb{N}$, $\convergenceorder>q$.
It holds 
$$ \left\lvert \frac{(c\stepsize)^{q+1}}{\lceil \lvert c \rvert \cdot \mrfactor{} \rceil^{q}} \right\rvert = \mathcal{O}(\stepsize^{\convergenceorder+1}), $$
provided that the multirate factor satisfies
$\stepsize^{(q-p)/q} \leq \mrfactor{}.$
}
\begin{proof}
    By Lemma \ref{lemma:mr-factor_optimality} and the lower bound on $\mrfactor{}$, we obtain 
    \begin{equation*}
    \begin{aligned}
        \left\lvert \frac{(c\stepsize)^{q+1}}{\lceil \lvert c \rvert \cdot \mrfactor{} \rceil^{q}}  \right\rvert
    &= 
    \lvert c \rvert \cdot \left(\frac{\lvert c \rvert \cdot \stepsize}{\lceil \lvert c \rvert \cdot \mrfactor{} \rceil}\right)^q  \stepsize 
    \leq 
    \lvert c \rvert \cdot \left(\frac{\stepsize}{\mrfactor{}}\right)^q \stepsize
    =
    \lvert c \rvert \cdot \frac{\stepsize^{q+1}}{\mrfactor{}^q}  
    \\
    &\leq 
    \lvert c \rvert \cdot (\stepsize^{(\convergenceorder-q)/q})^q \cdot \stepsize^{q+1} 
    = 
    \lvert c \rvert \cdot \stepsize^{\convergenceorder + 1} 
    = 
    \mathcal{O}(\stepsize^{\convergenceorder + 1}). & \qedhere
    \end{aligned}
    \end{equation*}
\end{proof}\normalfont

Roughly speaking, if we replace for all nodes $\node{v} \in \nodeset \setminus \{\node{r}\}$ the factor $\mrfactor{\node{v}}$ with 
\begin{equation}\label{eq:step-size_dependent_mrfactors}
    \mrfactor{\node{v}}(\stepsize) = \stepsize^{(\convergenceorder^{\{\node{v}\}} - \convergenceorder^{\{\node{r}\}})/\convergenceorder^{\{\node{v}\}}},
\end{equation}
where $\convergenceorder^{\{\node{r}\}}$ is the order of the method in root node $\node{r}$, the multirate hierarchical splitting method scales with order $\convergenceorder^{\{\node{r}\}}$. 
To maintain a choice of multirate factors independent of $\stepsize$ and $\convergenceorder$, we arrive at the following definition of computational order.
\definition[Computational order]{
A multirate hierarchical splitting method $\SplittingMethod{\stepsize}{}$ exhibits a computational order of $\convergenceorder$ for step sizes $\stepsize \geq \stepsize_{\mathtt{min}}$, provided that $\SplittingMethod{\stepsize}{\node{r}}$ has order $\convergenceorder^{\{\node{r}\}} = \convergenceorder$; and for each node $\node{v} \in \nodeset \setminus \{\node{r}\}$ with $\convergenceorder^{\{\node{v}\}} < \convergenceorder^{\{\node{r}\}}$, $\mrfactor{\node{v}} \geq \stepsize_{\mathtt{min}}^{(\convergenceorder^{\{\node{v}\}} - \convergenceorder^{\{\node{r}\}})/\convergenceorder^{\{\node{v}\}}}$.
}\normalfont\medskip

A multirate hierarchical splitting method of computational order $\convergenceorder$ is expected to be more accurate for step sizes $\stepsize \geq \stepsize_{\mathtt{min}}$ compared to its version employing multirate factors \eqref{eq:step-size_dependent_mrfactors}, which scales with order $\convergenceorder$.
The concept of computational order has already been applied in molecular dynamics simulations~\cite{shcherbakov2017schwinger}.\footnote{In \cite{shcherbakov2017schwinger}, the notion of computational order was referred to as \emph{effective order}. In the sense of \cite{BUTCHER1998179}, however, a method has effective order $\convergenceorder$ if a pre- and post-processing turns it into a method of order $\convergenceorder$, as exploited by processed splitting methods \cite{blanes2006processing}. To avoid any ambiguity with this established concept, and since the notion concerns the superior accuracy observed in practical computations over a limited range of step sizes, we use the term \emph{computational order}.}

\section{Design of hierarchical splitting methods}\label{sec:design}

In this section, we discuss possibilities for setting up the splitting tree for $\Npartitions$-split systems.
We start by demonstrating that state-of-the-art splitting methods for $\Npartitions$-split systems can be represented as hierarchical splitting methods.
Subsequently, in \Cref{sec:design_principles}, we will discuss principles for designing splitting trees, with particular focus on two important subclasses of dynamical systems: separable Hamiltonian systems and port-Hamiltonian systems of ordinary differential equations.

\subsection{Representation of state-of-the-art methods}
The state-of-the-art splitting methods for $\Npartitions$-split systems introduced in \Cref{sec:introduction} can all be represented as hierarchical splitting methods.
The Lie--Trotter splitting~\eqref{eq:Lie-Trotter} and the Strang splitting~\eqref{eq:Strang-splitting} correspond to the splitting tree from \Cref{fig:mr_HSM_splitting_tree} with multirate factor~$1$ on each edge, the respective splitting method for two-split systems in all inner nodes, and the exact flows in all leaf nodes.
Considering a base method $\SplittingMethod{\stepsize}{}$ of even order $\convergenceorder$, a composition scheme~\eqref{eq:composition_based_on_Strang} with weights $\gamma_i$ for $i=1,\ldots,s$ satisfying $\sum_i \gamma_i = 1$ and $\sum_i (\gamma_i)^{\convergenceorder + 1} = 0$ is convergent of order at least $\convergenceorder+1$.
If, in addition, the base method is self-adjoint and the weights are palindromic, $\gamma_i = \gamma_{s+1-i}$ for $i=1,\ldots,s$, the composition is itself self-adjoint and therefore convergent of order $\convergenceorder+2$.
One example is given by the \emph{triple-jump composition} \cite{yoshida1990construction}
\begin{equation}\label{eq:triple-jump_weights}
    \begin{aligned}
        \gamma_1 &= \gamma_3 = 1 / (2 - 2^{1/3}), & \gamma_2 &= 1-2\gamma_1.
    \end{aligned}
\end{equation}
When choosing the Strang splitting \eqref{eq:Strang_N2} as the base method, the composition scheme \eqref{eq:composition_based_on_Strang} with weights \eqref{eq:triple-jump_weights} is convergent of order $\convergenceorder = 4$ and reads
\begin{subequations}\label{eq:Yoshida_both-versions}
\begin{equation}\label{eq:Yoshida}
    \SplittingMethod{\stepsize}{} = \flow{\gamma_1\stepsize/2}{1} \circ \flow{\gamma_1\stepsize}{2} \circ \flow{\gamma_1\stepsize/2}{1} \circ \flow{\gamma_2\stepsize/2}{1} \circ \flow{\gamma_2\stepsize}{2} \circ \flow{\gamma_2\stepsize/2}{1} \circ \flow{\gamma_1\stepsize/2}{1} \circ \flow{\gamma_1\stepsize}{2} \circ \flow{\gamma_1\stepsize/2}{1}.
\end{equation}
If the flows of the subsystems are computed exactly, \eqref{eq:Yoshida} is equivalent to
\begin{equation}\label{eq:Yoshida_compressed}
    \SplittingMethod{\stepsize}{} = \flow{\gamma_1\stepsize/2}{1} \circ \flow{\gamma_1\stepsize}{2} \circ \flow{(\gamma_1+\gamma_2)\stepsize/2}{1} \circ \flow{\gamma_2\stepsize}{2} \circ \flow{(\gamma_1+\gamma_2)\stepsize/2}{1} \circ \flow{\gamma_1\stepsize}{2} \circ \flow{\gamma_1\stepsize/2}{1}.
\end{equation}
\end{subequations}
Hierarchical splitting methods $\SplittingMethod{\stepsize}{}$ of even order $\convergenceorder$ can also be employed as the base method for the composition~\eqref{eq:composition_based_on_Strang} to derive higher-order splitting methods for $\Npartitions$-split systems.
The resulting composition scheme is then once again a hierarchical splitting method.
\lemma{
Let $\SplittingMethod{\stepsize}{}$ be a hierarchical splitting method.
Then any composition method \eqref{eq:composition_based_on_Strang} using $\SplittingMethod{\stepsize}{}$ as the base method is again a hierarchical splitting method, whose splitting tree differs only in the root node, where $\SplittingMethod{\stepsize}{\node{r}}$ is replaced by its composition with the same weights.
}
\begin{proof}
    The composition of the full hierarchical splitting method reads 
    $$ \bm{\Psi}_{\stepsize} = \SplittingMethod{\gamma_s \stepsize}{} \circ \SplittingMethod{\gamma_{s-1} \stepsize}{} \circ \ldots \circ \SplittingMethod{\gamma_{1} \stepsize}{}. $$
    By construction, any application of the base method starts with an application of the method $\SplittingMethod{\gamma_i \stepsize}{\node{r}}$ assigned to the root node.
    Consequently, the initial step of the iterative splitting yields
    $$\SplittingMethod{\gamma_s \stepsize}{\node{r}} \circ \SplittingMethod{\gamma_{s-1} \stepsize}{\node{r}} \circ \ldots \circ \SplittingMethod{\gamma_{1} \stepsize}{\node{r}},$$
    which is equivalent to directly applying the composition technique to the method $\SplittingMethod{\stepsize}{\node{r}}$ assigned to the root node.
\end{proof}\normalfont

When considering compositions of self-adjoint splitting methods, the last flow of the $i$-th application and the first flow of the $(i+1)$-th application of the base scheme coincide. 
As it holds $\flow{\gamma_{i+1} \stepsize}{\node{v}} \circ \flow{\gamma_i \stepsize}{\node{v}} = \flow{(\gamma_i + \gamma_{i+1})\stepsize}{\node{v}}$ for the exact flows, these subsequent flow evaluations are frequently combined to reduce computational cost, see for example \eqref{eq:Yoshida_compressed}.
In a hierarchical splitting method, however, we do in general not assign the exact flows to the child nodes of the root node.
Consequently, combining subsequent flows results in a different method that may exhibit a lower order of convergence.
\example
{\label{ex:composition_HSM}
Consider the hierarchical splitting method \eqref{eq:example_HSM} and its splitting tree from \Cref{fig:N5_splitting-tree}. 
When replacing the Strang splitting in the root with the composition method \eqref{eq:Yoshida}, \Cref{theorem:convergence_order} states that all error terms $\mathcal{O}_3^{\{\node{i}\}}$ for $\node{i} \in \innernodes \setminus \{\node{r}\}$ are scaled by 
\begin{equation*} 
\begin{aligned}
&2 \left( \frac{\gamma_1}{2}\right)^3 + 2 \left( \frac{\gamma_2}{2}\right)^3 + 2 \left( \frac{\gamma_3}{2}\right)^3 = \frac{1}{4} \big( \gamma_1^3 + \gamma_2^3 + \gamma_3^3 \big) = 0, &
&\text{and} & 
&\gamma_1^3 + \gamma_2^3 + \gamma_3^3 = 0,
\end{aligned}
\end{equation*}
in the left and right branch, respectively. 
Hence, the third-order error terms all vanish.
Due to the inherent symmetry of the method, the hierarchical splitting method is convergent of order $\convergenceorder = 4$.
Alternatively, we may replace the Strang splitting in the root node with the composition method \eqref{eq:Yoshida_compressed}.
The leading error terms in the right branch of the splitting tree are scaled by $\gamma_1^3 + \gamma_2^3 + \gamma_3^3 = 0$ and thus vanish. However, the leading error terms in the left branch are scaled by
$$ \left(\frac{\gamma_1}{2}\right)^3 + \left(\frac{\gamma_1 + \gamma_2}{2}\right)^3 + \left(\frac{\gamma_2 + \gamma_3}{2}\right)^3 + \left(\frac{\gamma_3}{2}\right)^3  = \frac{3}{8} \left( \gamma_1^2 \gamma_2 + \gamma_1 \gamma_2^2 + \gamma_2^2 \gamma_3 + \gamma_2 \gamma_3^2 \right) \neq 0.$$
As a result, the overall hierarchical splitting method is only convergent of order $\convergenceorder = 2$.
}\normalfont\medskip

Note that Example \ref{ex:composition_HSM} includes a hierarchical splitting method of order $\convergenceorder$, wherein splitting methods of a lower convergence order $q < \convergenceorder$ are assigned to nodes of the splitting tree.

\lemma{\label[lemma]{lemma:composition_as_HSM}
Composition methods of the form \eqref{eq:composition_based_on_LT} using a first-order hierarchical splitting method $\bm{\chi}_{\stepsize}$ as the underlying base scheme can be represented as a hierarchical splitting method.}
\begin{proof}
Consider the coefficients $\alpha_i$ and $\beta_i$ from the composition \eqref{eq:composition_based_on_LT} and define $c \coloneqq \sum_{i=1}^s \alpha_i$, which we assume to satisfy $c \notin \{0,1\}$. 
Provided that the composition is consistent, it holds $\sum_{i=1}^s \alpha_i + \beta_i = 1$ and hence $\sum_{i=1}^s \beta_i = 1-c$.
Starting from the overall system \eqref{eq:split-system}, we consider the initial splitting 
$$\rhs{}(\x) = c \rhs{}(\x) + (1-c) \rhs{}(\x),$$
with exact flows $\flow{\stepsize}{1} = \flow{c\stepsize}{}$ and $\flow{\stepsize}{2} = \flow{(1-c)\stepsize}{}$, where $\flow{}{}$ denotes the exact flow of \eqref{eq:split-system}.
We define the scaled coefficients $\alpha_i^\prime \coloneqq \alpha_i/c$ and $\beta_i^\prime \coloneqq \beta_i/(1-c)$, so that $\flow{\alpha_i^\prime \stepsize}{1} = \flow{\alpha_i \stepsize}{}$ and $\flow{\beta_i^\prime \stepsize}{2} = \flow{\beta_i \stepsize}{}$.
To the root node, we assign the splitting method 
\begin{equation}\label{eq:LT_composition_as_HSM}
    \SplittingMethod{\stepsize}{\node{r}} = \flow{\beta_s^\prime \stepsize}{2} \circ \flow{\alpha_s^\prime \stepsize}{1} \circ \ldots \circ \flow{\beta_1^\prime \stepsize}{2} \circ \flow{\alpha_1^\prime \stepsize}{1}.
\end{equation}
Assigning the adjoint $\bm{\chi}_{\stepsize}^*$ of the first-order hierarchical splitting method to the left child node $\{1\}$ and the base method $\bm{\chi}_{\stepsize}$ to the right child node $\{2\}$ replaces $\flow{\alpha_i^\prime \stepsize}{1}$ by $\bm{\chi}^*_{\alpha_i \stepsize}$ and $\flow{\beta_i^\prime \stepsize}{2}$ by $\bm{\chi}_{\beta_i \stepsize}$.
Hence, \eqref{eq:LT_composition_as_HSM} is exactly the composition method \eqref{eq:composition_based_on_LT}.
\end{proof}\normalfont

\remark{\label{remark:composition_and_multirate}
The use of the reweighting \eqref{eq:reweighted_mr-factor} presents challenges for composition techniques.  
In contrast to (singlerate) hierarchical splitting methods, a composition \eqref{eq:composition_based_on_Strang} based on a multirate hierarchical splitting method is not equivalent to applying the composition \eqref{eq:composition_based_on_Strang} to the method $\SplittingMethod{\stepsize}{\node{r}}$ assigned to the root node.
In a composition of a multirate hierarchical splitting method, each application of the base scheme uses the same reweighted multirate factors, whereas the weights $\gamma_i$ enter the reweighting when applying the composition to the root node. 
In the latter case, the effective base scheme changes in every stage, such that the order conditions of the composition no longer remove the leading error terms of the subtrees below a reweighted edge.
In this case, the order of the composition is retained only if the methods assigned to these subtrees are already of the target order of the composition.
Similar issues arise when considering compositions of the form \eqref{eq:composition_based_on_LT}.}

\subsection{Design principles}\label{sec:design_principles}
In general, decomposing the vector field of a dynamical system $\x^\prime = \rhs{}(\x)$ into an $\Npartitions$-split system \eqref{eq:split-system} constitutes an open challenge. 
Even if a decomposition of $\rhs{}$ into $\Npartitions$ components is given, it is well-known that the order of the operators inside a splitting method \eqref{eq:splitting_method} can significantly impact the accuracy and numerical stability, see e.g.~\cite{McLachlan2021Tuning}.
In this work, we do not aim at addressing these open challenges. 
Instead, we assume that, by physical or mathematical structure, a decomposition into $\Npartitions$ partitions is given and we have some knowledge about the properties and dynamics of the subsystems. 

\subsubsection{General guidelines}
In general, the number of sub-steps involving a numerical method $\SplittingMethod{\stepsize}{\node{v}}$ increases exponentially with the distance of the node $\node{v}$ from the root node.
If the system does not exhibit multirate potential that justifies certain subflows being computed more frequently, the exponentially growing number of sub-steps can become a significant drawback of the hierarchical construction, leading to an inefficient computational process. In such situations, it is advisable to maintain a small depth of the splitting tree by designing balanced binary trees.

For certain classes of dynamical systems, there exist commutator-based splitting approaches \cite{chin1997symplectic,blanes2023ModifiedPotentials,MONCH202525,Monch2026splitting} that enhance the accuracy compared to conventional splitting methods. 
The hierarchical splitting approach allows for assigning commutator-based splitting approaches to inner nodes of a splitting tree, provided that the approach is consistent with the Lie derivatives of the child nodes, cf.~\Cref{remark:commutators}.
In particular, commutator-based splitting methods typically have fewer stages compared to conventional splitting methods of the same order. 
Therefore, the use of commutator-based methods contributes to reducing the number of flow evaluations with increasing depth of the splitting tree.
When considering dynamical systems exhibiting multirate potential, the design of splitting trees can be fundamentally different. 
Since we are interested in computing fast sub-flows more frequently, the general guideline of designing balanced trees to keep the depth small may no longer apply. 
Instead, maximizing the depth of the splitting tree for a given number of partitions may be of particular interest, since it aligns the tree with the multirate potential of the system.

\remark{\label{remark:superiority}
Since the Lie--Trotter \eqref{eq:Lie-Trotter} and Strang \eqref{eq:Strang-splitting} splittings, as well as their compositions \eqref{eq:composition_based_on_LT} and \eqref{eq:composition_based_on_Strang}, can be represented as hierarchical splitting methods, the proposed framework naturally includes general-purpose methods for $\Npartitions$-split systems while simultaneously extending them by multiple time stepping techniques.
Hierarchical splitting is particularly advantageous whenever a) the system exhibits multirate potential, in which case the depth of the splitting tree can be tailored to the time-scale separation, b) commutator-based methods are applicable to the components assigned to the child nodes, reducing the number of stages relative to conventional splitting methods of the same order, or c) the subsystems admit dedicated methods of higher order or with structure-preserving properties, which a composition of a low-order base scheme cannot provide.
}\normalfont\medskip

\subsubsection{Separable Hamiltonian systems}
In molecular dynamics simulations, for instance in lattice quantum chromodynamics \cite{SEXTON1992665,URBACH200687}, Hamiltonian systems with fine structure 
\begin{equation}\label{eq:separable_H}
    \Hamiltonian(\p,\q) = \Kinetic(\p) + \Potential(\q) = \tfrac{1}{2} \p^\top \bm{M}^{-1} \p + \sum_{m=1}^{\Npartitions -1} \Potential^{\{m\}}(\q), \quad \bm{M} = \bm{M}^\top \succ 0,
\end{equation}
often arise.
Typically, the decomposition of the potential energy $\Potential$ gives rise to a force hierarchy, where the components $\Potential^{\{m\}}$ are ordered such that their respective forces are characterized by increasing contribution to the dynamics of the solution and simultaneously decreasing evaluation cost. 
This introduces a source of multirate potential, which motivates the use of multiple time stepping techniques. 
Besides the impulse method~\eqref{eq:impulse_method}, nested integration \cite{tuckerman_reversible_1992,SEXTON1992665} constitutes a class of multiple time stepping schemes frequently employed in molecular dynamics simulations. 
By decomposing \eqref{eq:separable_H} into its $\Npartitions$ summands and considering $\flow{\stepsize}{m}$ as the exact flows of the subsystems $\p^\prime = - \nabla \Potential^{\{m\}}(\q)$, $\q^\prime = \bm{0}$, and $\flow{\stepsize}{\Npartitions}$ the exact flow of the subsystem $\p^\prime = \bm{0}$, $\q^\prime = \bm{M}^{-1} \p$, nested integration schemes are defined as 
\begin{equation}\label{eq:nested_integration}
\begin{aligned} 
    \SplittingMethod{\stepsize}{} &= \flow{b_{s^{\{1\}}}^{\{1\}} \stepsize}{1} \circ \Big( \SplittingMethod{a_{s^{\{1\}}}^{\{1\}} \stepsize/\mrfactor{1}}{2} \Big)^{\mrfactor{1}} \circ \ldots \circ \flow{b_1^{\{1\}} \stepsize}{1} \circ \Big(\SplittingMethod{a_1^{\{1\}} \stepsize/\mrfactor{1}}{2}\Big)^{\mrfactor{1}}, \\
    \SplittingMethod{\tau}{2} &= \flow{b_{s^{\{2\}}}^{\{2\}} \tau}{2} \circ \Big( \SplittingMethod{a_{s^{\{2\}}}^{\{2\}} \tau/\mrfactor{2}}{3} \Big)^{\mrfactor{2}} \circ \ldots \circ \flow{b_1^{\{2\}} \tau}{2} \circ \Big(\SplittingMethod{a_1^{\{2\}} \tau/\mrfactor{2}}{3}\Big)^{\mrfactor{2}}, \\
    &\;\;\vdots \\
    \SplittingMethod{\tau}{\Npartitions - 1} &= \flow{b_{s^{\{\Npartitions - 1\}}}^{\{\Npartitions - 1\}} \tau}{\Npartitions-1} \circ \flow{a_{s^{\{\Npartitions - 1\}}}^{\{\Npartitions-1\}} \tau}{\Npartitions} \circ \ldots \circ \flow{b_1^{\{\Npartitions - 1\}} \tau}{\Npartitions-1} \circ \flow{a_1^{\{\Npartitions-1\}} \tau}{\Npartitions}.
\end{aligned}
\end{equation}
This exactly corresponds to the hierarchical splitting approach with the splitting tree from \Cref{fig:mr_HSM_splitting_tree}, where the labels $m=1,\ldots,\Npartitions-1$ correspond to the components $\Potential^{\{m\}}$, and label $\Npartitions$ to the kinetic energy $\Kinetic$.
For numerical integration of separable Hamiltonian systems with quadratic kinetic energy $\Kinetic$, force-gradient integrators \cite{chin1997symplectic,omelyan2003symplectic}, also known as splitting methods based on modified potentials \cite{blanes2023ModifiedPotentials}, incorporate the force-gradient term 
\begin{equation}\label{eq:fg_term}
    \big[ \LieDerivative{\Potential} , [ \LieDerivative{\Kinetic} , \LieDerivative{\Potential} ] \big] = 2 \nabla \Potential(\q) \bm{M}^{-1} \nabla^2 \Potential(\q) \frac{\partial}{\partial \p}
\end{equation}
into the computational process to enhance the accuracy. 
In the inner nodes of the splitting tree, the corresponding commutator for $\SplittingMethod{}{m}$ generalizes to 
\begin{equation}\label{eq:generalized_fg-term}
    \Big[ \LieDerivative{\Potential^{\{m\}}} , \Big[ \LieDerivative{\Kinetic} + {\textstyle \sum_{k=m+1}^{\Npartitions - 1}} \LieDerivative{\Potential^{\{k\}}} , \LieDerivative{\Potential^{\{m\}}} \Big] \Big] .
\end{equation}
Since the vector fields $\LieDerivative{\Potential^{\{m\}}}$ for $m=1,\ldots,\Npartitions-1$ act only on $\p$ and depend only on $\q$, they pairwise commute, $[\LieDerivative{\Potential^{\{m\}}},\LieDerivative{\Potential^{\{k\}}}] \equiv 0$.
Consequently, \eqref{eq:generalized_fg-term} simplifies to 
\begin{equation*}
    \big[ \LieDerivative{\Potential^{\{m\}}} , [ \LieDerivative{\Kinetic} , \LieDerivative{\Potential^{\{m\}}} ] \big],
\end{equation*}
such that the force-gradient approach is compatible with nested integration and the method $\SplittingMethod{}{m}$ only requires information of the local component $\Potential^{\{m\}}$.
Therefore, it is possible to assign force-gradient integrators to inner nodes of the splitting tree in \Cref{fig:mr_HSM_splitting_tree}, thereby obtaining multiple time stepping extensions of force-gradient integrators.

Embedding nested integration schemes into the hierarchical splitting framework thus yields a refined error analysis, reduces the computational overhead via the reweighting, and makes force-gradient integrators available in the inner nodes. 
These enhancements are illustrated numerically in \Cref{sec:FPU}, for a modified Fermi--Pasta--Ulam (FPU) problem.

\subsubsection{Port-Hamiltonian systems}\label{sec:design_PHS}
Port-Hamiltonian systems of ordinary differential equations (pH-ODEs) often arise in the modeling of dynamical systems~\cite{VanDerSchaft2014PHS} and take the form 
\begin{equation}\label{eq:PHS}
\begin{aligned}
    \x^\prime &= \big[\bm{J}(\x) - \bm{R}(\x) \big] \nabla \mathcal{H}(\x) + \bm{B}(\x) \bm{u}(t), &
    \bm{y} &= \bm{B}(\x)^\top \nabla \mathcal{H}(\x),
\end{aligned}
\end{equation}
with point-wise skew-symmetric interconnection matrix $\bm{J}(\x) \in \mathbb{R}^{n \times n}$ ($\bm{J} = -\bm{J}^\top$), point-wise symmetric and positive semi-definite dissipation matrix $\bm{R}(\x) \in \mathbb{R}^{n \times n}$ ($\bm{R} = \bm{R}^\top \succeq \bm{0}$), port matrix $\bm{B}(\x) \in \mathbb{R}^{n \times m}$, input $\bm{u}(t) \in \mathbb{R}^m$, and corresponding output $\bm{y}(t) \in \mathbb{R}^m$.  
An important feature of port-Hamiltonian systems is the energy exchange through their ports. 
Assuming sufficient regularity of the input, the solution satisfies the power balance 
\begin{equation}\label{eq:PBE}
    \tfrac{\mathrm{d}}{\mathrm{d} t} \mathcal{H}(\x(t)) = -\nabla \mathcal{H}(\x(t))^\top \bm{R}(\x(t)) \nabla \mathcal{H}(\x(t)) + \bm{y}^\top \bm{u}(t) \leq \bm{y}^\top \bm{u}(t).
\end{equation}
At the discrete level, this structure should be reflected by employing energy-consistent numerical integration schemes.
A splitting method is energy-consistent, provided that the subsystems are pH-ODEs with respect to the Hamiltonian $\Hamiltonian$ of the original system, all subsystems involving dissipation are only integrated forward in time, and the subflows are computed exactly or numerically via an energy-consistent method \cite{moench2026decomposition}. 
Energy-consistent numerical integration schemes for pH-ODEs include discrete gradient methods \cite{KINON202645,karsai2026DGM,Lautwein2026DGM} and Petrov--Galerkin methods \cite{EggerHabrichShashkov2021,giesselmann2024energyconsistentpetrovgalerkin}.

\noindent\begin{minipage}[c]{.735\textwidth}
Without imposing any further assumptions on the structure matrices $\bm{J}, \bm{R}$, and $\bm{B}$, we have a natural decomposition of \eqref{eq:PHS} into three port-Hamiltonian subsystems: $\x^\prime = \bm{J}(\x) \nabla \Hamiltonian(\x)$, $\x^\prime = - \bm{R}(\x) \nabla \Hamiltonian(\x)$, and $\x^\prime = \bm{B}(\x) \bm{u}(t)$.
In the spirit of hierarchical splitting, decomposing the port-Hamiltonian system amounts to a combination of the energy-associated and the port-based decomposition, cf.~\cite{moench2026decomposition}.
For linear port-Hamiltonian systems with constant structure matrices and quadratic Hamiltonian, applying the port-based decomposition first, followed by the energy-associated decomposition, enables the construction of fourth-order energy-consistent splitting methods by using (commutator-based) splitting methods for two-split systems from~\cite{MONCH202525}.
The resulting splitting tree is illustrated in~\Cref{fig:phs_splitting_tree}.
\end{minipage}
\hfill 
\begin{minipage}[c]{.215\textwidth}
\vfill
\begin{figure}[H]
    \centering
    \begin{tikzpicture}[level distance=1.25cm,
          level 1/.style={sibling distance=1.2cm},
          level 2/.style={sibling distance=1.2cm},
          level 3/.style={sibling distance=1.2cm}]
          \node {$\{\bm{J},\bm{R},\bm{B}\}$}
            child {
            node {$\{\bm{J},\bm{R}\}$}
              child {node {$\{\bm{J}\}$}
              }
              child {node {$\{\bm{R}\}$}
              }
            }
            child {node {$\{\bm{B}\}$
            }
            };
        \end{tikzpicture}
    \caption{Splitting tree for pH-ODEs.}
    \label{fig:phs_splitting_tree}
\end{figure}
\vfill
\end{minipage}

\vspace*{-3ex}

\remark{\label{remark:coupled_pHS}
An important structural benefit of port-Hamiltonian systems is that appropriate coupling of port-Hamiltonian systems yields a monolithic port-Hamiltonian system \cite{ehrhardt2025coupling}. 
Such systems typically arise in multi-physics applications, where subsystems frequently evolve on distinct time scales and require multiple time stepping. 
Given a coupled pH-ODE, the block structure of the system matrices $\bm{J}$, $\bm{R}$, and $\bm{B}$ offers additional possibilities for structure-preserving splitting.
A comprehensive discussion of (hierarchical) splitting approaches for coupled port-Hamiltonian systems is beyond the scope of the present work.
In \Cref{sec:MSD}, however, we will consider a coupled pH-ODE with multirate potential that highlights the structural benefits of the hierarchical splitting approach in designing higher-order energy-consistent splitting methods including multiple time stepping.
}

\section{Numerical results}\label{sec:NumericalResults}
In this section, we examine the performance of hierarchical splitting methods when applied to different test problems.
The numerical simulations have been performed in Python and the corresponding Jupyter notebooks are available on GitHub, \url{https://github.com/KevinSchaefers/HierarchicalSplittingMethods}.

\subsection{Fermi--Pasta--Ulam problem}\label{sec:FPU}
To illustrate the enhancements for nested integration schemes provided by the hierarchical splitting approach, we consider the modified Fermi--Pasta--Ulam (FPU) problem introduced by \cite{GALGANI1992334}.
It describes a chain of $2m$ mass points that are connected alternately by stiff linear springs and soft nonlinear springs, where we consider $m = 3$ throughout this section.
In the variables $\q = (\q^{\{\mathfrak{s}\}},\q^{\{\mathfrak{f}\}})$ and $\p = (\p^{\{\mathfrak{s}\}},\p^{\{\mathfrak{f}\}})$ with $\q^{\{\mathfrak{s}\}},\q^{\{\mathfrak{f}\}},\p^{\{\mathfrak{s}\}},\p^{\{\mathfrak{f}\}} \in \mathbb{R}^m$, which describe the motion of the centers of mass of the stiff springs and their elongations, the system is modeled via the separable Hamiltonian
\begin{equation}\label{eq:FPU_Hamiltonian}
\begin{aligned}
    \Hamiltonian(\p,\q) &= \Kinetic(\p) + \Potential^{\{\mathfrak{f}\}}(\q) + \Potential^{\{\mathfrak{s}\}}(\q), \quad \Kinetic(\p) = \tfrac12 \lVert \p \rVert_2^2, \quad
        \Potential^{\{\mathfrak{f}\}}(\q) = \tfrac{\omega^2}{2} \lVert \q^{\{\mathfrak{f}\}} \rVert_2^2, \\
        \Potential^{\{\mathfrak{s}\}}(\q) &= \tfrac{1}{4} \big( (q^{\{\mathfrak{s}\}}_{1} - q^{\{\mathfrak{f}\}}_{1})^4 + (q^{\{\mathfrak{s}\}}_{m} + q^{\{\mathfrak{f}\}}_{m})^4 + {\textstyle\sum_{i=1}^{m-1}} \big(q^{\{\mathfrak{s}\}}_{i+1} - q^{\{\mathfrak{f}\}}_{i+1} - q^{\{\mathfrak{s}\}}_{i} - q^{\{\mathfrak{f}\}}_{i}\big)^4 \big).
\end{aligned}
\end{equation}
The stiff springs give rise to the fast potential $\Potential^{\{\mathfrak{f}\}}$ with the frequency $\omega = 500$, whereas the soft nonlinear springs form the slow potential $\Potential^{\{\mathfrak{s}\}}$.
As initial values, we choose $q^{\{\mathfrak{s}\}}_{1}(0) = 1$, $q^{\{\mathfrak{f}\}}_{1}(0) = \omega^{-1}$, $p^{\{\mathfrak{s}\}}_{1}(0) = p^{\{\mathfrak{f}\}}_{1}(0) = 1$, and zero for all other components.
We conduct numerical simulations on the time interval $t \in [0,10]$.

For designing hierarchical splitting methods, we consider the splitting tree of nested integration schemes (cf.~\Cref{fig:mr_HSM_splitting_tree}) for the three partitions $\Kinetic$, $\Potential^{\{\mathfrak{f}\}}$, and $\Potential^{\{\mathfrak{s}\}}$: the root node splits the slow potential $\Potential^{\{\mathfrak{s}\}}$ from the fast subsystem $\Kinetic + \Potential^{\{\mathfrak{f}\}}$, which is split further into $\Kinetic$ and $\Potential^{\{\mathfrak{f}\}}$, and the multirate factor $\mrfactor{}$ is assigned to the edge between root node and node $\Kinetic + \Potential^{\{\mathfrak{f}\}}$.
To all leaf nodes, the exact flows of the subsystems are assigned.
The FPU problem is a popular test problem for multiple time stepping methods.
As we aim at comparing nested integration schemes, we do not solve the fast subsystem $\Kinetic + \Potential^{\{\mathfrak{f}\}}$ exactly, although this would be possible due to its linearity.
We consider four hierarchical splitting methods that differ in the splitting methods assigned to the inner nodes.
The method \texttt{HOMF4} uses the fourth-order method \eqref{eq:OMF4} in both nodes.
The method \texttt{COMP4} uses \eqref{eq:OMF4} in the root node and the Strang splitting \eqref{eq:Strang_N2} in the fast node, such that it is only convergent of order two, but exhibits a computational order of four for sufficiently large $\mrfactor{}$.
The method \texttt{HOMF6} uses the sixth-order method \eqref{eq:OMF6} in both nodes.
The method \texttt{HFG4} uses the fourth-order force-gradient method \eqref{eq:ABACABACABA} in both nodes, where the force-gradient terms are evaluated exactly by Hessian--vector products.
For comparison, we also perform numerical tests with the state-of-the-art composition schemes \texttt{BMs6p4} \eqref{eq:BMs6p4} and \texttt{BMs10p6} \eqref{eq:BMs10p6} of the form \eqref{eq:composition_based_on_LT} with
\begin{equation}\label{eq:FPU_LT-splitting}
    \bm{\chi}_{\stepsize} = \big( \flow{\stepsize/\mrfactor{}}{\Kinetic} \circ \flow{\stepsize/\mrfactor{}}{\Potential^{\{\mathfrak{f}\}}}\big)^{\mrfactor{}} \circ  \flow{\stepsize}{\Potential^{\{\mathfrak{s}\}}}
\end{equation}
as the underlying base scheme.
Note that, as per Lemma~\ref{lemma:composition_as_HSM}, both \texttt{BMs6p4} and \texttt{BMs10p6} can be represented as hierarchical splitting methods.
The error of a numerical solution $\x_0,\ldots,\x_{N_t}$ with $\x = (\p,\q)$ and $\stepsize = 10/N_t$ is measured in the $L^2$-norm $\big( \stepsize \sum_{n=0}^{N_t} \lVert \x_n - \x(t_n) \rVert_2^2 \big)^{1/2}$.

\begin{figure}[htb!]
    \centering

    \begin{subfigure}{0.475\textwidth}
        \centering
        \includegraphics{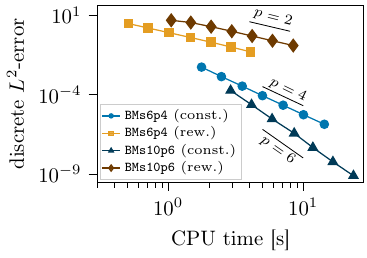}
        \vspace*{-2ex}
        \caption{\footnotesize Composition schemes \texttt{BMs6p4} and \texttt{BMs10p6}. Comparison for $\mrfactor{} = 10$ when applying the reweighting \eqref{eq:reweighted_mr-factor} vs.~no reweighting.}
        \label{fig:FPU_BM_results}
    \end{subfigure}
    \hfill
    \begin{subfigure}{0.475\textwidth}
        \centering
        \includegraphics{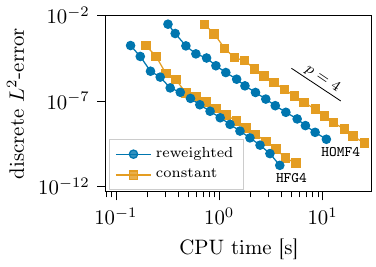}
        \vspace*{-2ex}
        \caption{\footnotesize \texttt{HOMF4} and \texttt{HFG4} methods. Comparison of the variants with reweighting ($\mrfactor{} = \mathrm{30}$ and $\mrfactor{} = 15$) and without reweighting ($\mrfactor{} = 17$ and $\mrfactor{} = 4$).}
        \label{fig:FPU_HOMF4_results}
    \end{subfigure}
    \\
    \begin{subfigure}{0.475\textwidth}
        \centering
        \includegraphics{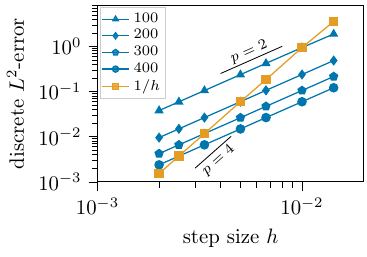}
        \vspace*{-2ex}
        \caption{\footnotesize \texttt{COMP4} method. Comparison for different multirate factors $\mrfactor{} \in \{100,200,300,400,1/\stepsize\}$ while applying the reweighting \eqref{eq:reweighted_mr-factor}.}
        \label{fig:FPU_COMP4_results}
    \end{subfigure}
    \hfill
    \begin{subfigure}{0.475\textwidth}
        \centering
        \includegraphics{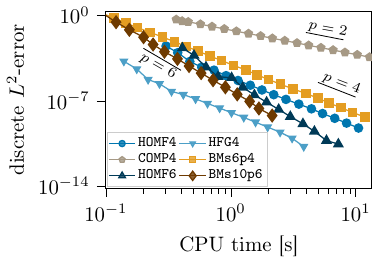}
        \vspace*{-2ex}
        \caption{\footnotesize Work-precision diagram for all methods with the multirate factors stated in the text. Discrete $L^2$-error vs.~CPU time (in seconds).}
        \label{fig:FPU_work-precision}
    \end{subfigure}

    \caption{\footnotesize Fermi--Pasta--Ulam problem. Numerical results for (a) the composition methods \texttt{BMs6p4} and \texttt{BMs10p6}, (b) the fourth-order methods \texttt{HOMF4} and \texttt{HFG4}, (c) the \texttt{COMP4} scheme, and (d) the computational efficiency of all methods.}
\end{figure}

\Cref{fig:FPU_BM_results} shows an order reduction to order two for both composition schemes\footnote{Due to the inherent symmetry of the methods, the order is two instead of the order of the underlying base scheme ($\convergenceorder = 1$).} when the reweighting \eqref{eq:reweighted_mr-factor} is applied to the multirate factors.
With the constant multirate factor $\mrfactor{} = 10$, the schemes attain their orders $\convergenceorder = 4$ and $\convergenceorder = 6$, respectively.
These observations are in line with \Cref{remark:composition_and_multirate}.

To quantify the savings from the reweighting \eqref{eq:reweighted_mr-factor}, we compare \texttt{HOMF4} and \texttt{HFG4} with reweighted multirate factors $\mrfactor{} = 30$ and $\mrfactor{} = 15$ against constant multirate factors chosen such that the sub-step of largest modulus is resolved with the same micro step size, leading to $\mrfactor{} = 17$ and $\mrfactor{} = 4$, respectively.
\Cref{fig:FPU_HOMF4_results} shows that the reweighting substantially enhances the efficiency of both methods: the reweighted variant requires roughly $50\,\%$ of the CPU time for \texttt{HOMF4} and roughly $75\,\%$ for \texttt{HFG4}, while, in contrast to the composition techniques, the fourth order is preserved in both cases.
Although reweighting the multirate factors results in an order reduction when applied to composition techniques, it contributes to minimizing computational overhead for hierarchical splitting methods that are not based on a composition of a lower-order base method.

\Cref{fig:FPU_COMP4_results} illustrates the concept of computational order for the \texttt{COMP4} method with the reweighting \eqref{eq:reweighted_mr-factor} and varying multirate factors.
A step size-dependent multirate factor $\mrfactor{} = 1/\stepsize$ yields a method that scales with order $\convergenceorder = 4$, whereas all fixed multirate factors yield order $\convergenceorder = 2$ while being more accurate for step sizes $\stepsize \geq 1/\mrfactor{}$.

\Cref{fig:FPU_work-precision} shows a work-precision diagram for the fourth-order methods \texttt{HOMF4} (reweighted $\mrfactor{} = 30$), \texttt{COMP4} (reweighted $M = 200$), \texttt{HFG4} (reweighted $\mrfactor{} = 15$), and \texttt{BMs6p4} (constant $\mrfactor{} = 1$), as well as the sixth-order methods \texttt{HOMF6} (reweighted $\mrfactor{} = 12$) and \texttt{BMs10p6} (constant $\mrfactor{} = 1$), where the multirate factors were tuned manually such that each variant achieves the best performance in the work-precision diagram.
Among the methods of order four, \texttt{HOMF4} requires roughly half the CPU time of \texttt{BMs6p4}, and exploiting the force-gradient approach is even more rewarding: \texttt{HFG4} needs about $11\,\%$ of the CPU time of \texttt{BMs6p4}.
The method \texttt{COMP4} is by far the least efficient scheme considered, since its computational order of four is only attained for step sizes $\stepsize \geq 1/\mrfactor{}$, whereas the relevant range of accuracies requires smaller step sizes.
Among the methods of order six, \texttt{HOMF6} requires about $20$ to $55\,\%$ more CPU time than \texttt{BMs10p6}.
\texttt{HFG4} is nevertheless the most efficient scheme over the entire range of accuracies considered.

\begin{figure}[htbp!]
    \centering
    \includegraphics{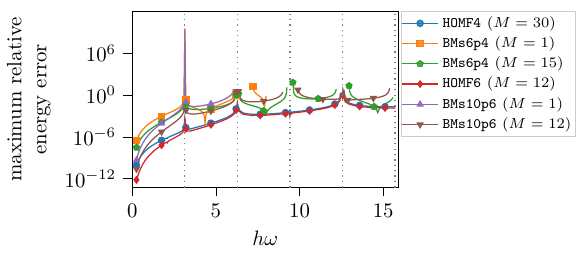}

    \caption{\footnotesize Fermi--Pasta--Ulam problem. Robustness with respect to the step size for the methods \texttt{HOMF4} (reweighted $\mrfactor{}=30$), \texttt{HOMF6} (reweighted $\mrfactor{}=12$), \texttt{BMs6p4} (constant $\mrfactor{}=1$ and $\mrfactor{}=15$), and \texttt{BMs10p6} (constant $\mrfactor{}=1$ and $\mrfactor{}=12$). Maximum relative energy error for $\stepsize\omega \in [0.25,5\pi]$; the dotted lines mark the resonances $\stepsize\omega = k\pi$, $k = 1,\ldots,5$.
    Overflowing runs are not depicted.}
    \label{fig:FPU_stability}
\end{figure}

Finally, we investigate the robustness of \texttt{HOMF4}, \texttt{HOMF6}, and the two composition schemes with respect to the step size.
For each method, we fix the multirate factor from \Cref{fig:FPU_work-precision} and vary the step size such that $\stepsize\omega \in [0.25, 5\pi]$.
Since the micro step of the fast node grows with $\stepsize$ for fixed $\mrfactor{}$, we additionally consider \texttt{BMs6p4} and \texttt{BMs10p6} with the constant multirate factors $\mrfactor{} = 15$ and $\mrfactor{} = 12$, leading to methods requiring about as many fast force evaluations per macro step as \texttt{HOMF4} and \texttt{HOMF6}, respectively.
As a measure of robustness, we use the maximum relative energy error $\max_n \lvert \Hamiltonian(\x_n) - \Hamiltonian(\x_0) \rvert / \lvert \Hamiltonian(\x_0) \rvert$.
It is well known that multiple time stepping methods may suffer from resonances at $\stepsize\omega = k\pi$, $k \in \mathbb{N}$ \cite{BIESIADECKI1993}, as the exact flow of the fast subsystem over one step then acts as $\pm \mathrm{Id}$ on the fast components.
\Cref{fig:FPU_stability} shows that all methods are affected by these resonances, but to a different extent.
The composition schemes with $\mrfactor{} = 1$ overflow for the majority of the step sizes considered, namely for $55\,\%$ and $62\,\%$, whereas the hierarchical splitting methods do not overflow at any step size and keep the maximum relative energy error below $5$ over the entire range.
Increasing the multirate factor of the composition schemes to the cost-equivalent values $\mrfactor{} = 15$ and $\mrfactor{} = 12$ reduces the fraction of overflowing step sizes to $8.3\,\%$ and $8.6\,\%$, but both schemes still overflow near the resonances $\stepsize\omega = k\pi$ with $k \geq 2$.
These observations can be explained by the fact that each composition stage integrates the system with the first-order base scheme~\eqref{eq:FPU_LT-splitting}, whereas the hierarchical splitting methods assign a method of order four and six, respectively, to the fast node and thus remain robust at comparable cost.

The FPU problem satisfies all three conditions of \Cref{remark:superiority}, and each is reflected in the results: the multirate potential is exploited through the reweighted multirate factors, the applicability of commutator-based methods makes \texttt{HFG4} the most efficient scheme considered, and assigning a higher-order method to the fast node is what keeps the hierarchical splitting methods robust across the resonances.

\subsection{Coupled mass-spring-damper chains}\label{sec:MSD}
As a second example, we consider a modification of the coupled mass-spring-damper (MSD) chains from \cite{bartel2023dynamic} in which the linear springs are replaced by nonlinear ones.
The system consists of two nonlinear MSD chains with $N_1 = 5$ and $N_2 = 45$ masses, whose first masses are coupled by a linear spring.
Let $p_{i,j}$ and $q_{i,j}$, $j=1,\ldots,N_i$, denote the momenta and the spring deflections of chain $i \in \{1,2\}$, and let $q_c$ denote the deflection of the coupling spring. 
Collecting these quantities in $\x_1 \coloneqq (p_{1,1},q_{1,1},\ldots,p_{1,N_1},q_{1,N_1},q_c)^\top \in \mathbb{R}^{2N_1 + 1}$ and $\x_2 \coloneqq (p_{2,1},q_{2,1},\ldots,p_{2,N_2},q_{2,N_2})^\top \in \mathbb{R}^{2N_2}$, the state $\x = (\x_1,\x_2) \in \mathbb{R}^n$, $n=2(N_1 + N_2) + 1 = 101$, satisfies the pH-ODE 
\begin{equation}\label{eq:coupled_MSD}
    \x^\prime = \big( \bm{J} - \bm{R} \big) \nabla \Hamiltonian(\x), \qquad \bm{J} = \begin{pmatrix} \bm{J}_1 & \bm{C} \\ -\bm{C}^\top & \bm{J}_2 \end{pmatrix}, \qquad \bm{R} = \operatorname{blkdiag}(\bm{R}_1,\bm{R}_2).
\end{equation}
The matrices $\bm{J}_1,\bm{R}_1 \in \mathbb{R}^{(2N_1+1)\times(2N_1+1)}$, $\bm{J}_2,\bm{R}_2 \in \mathbb{R}^{2N_2\times 2N_2}$, and $\bm{C} \in \mathbb{R}^{(2N_1+1)\times 2N_2}$, as well as the Hamiltonian $\Hamiltonian$, are given by
\begin{align*}
    \bm{J}_1 &\coloneqq \begin{pmatrix}
        \operatorname{tridiag}(1,0,-1) & -\bm{e}_1 \\ \bm{e}_1^\top & 0
    \end{pmatrix}, & (\bm{R}_1)_{k,\ell} &\coloneqq \begin{cases} r_1, & k = \ell = 2j - 1, \; j = 1, \ldots, N_1, \\ 0, & \text{otherwise}, \end{cases} \\
    \bm{J}_2 &\coloneqq \operatorname{tridiag}(1,0,-1), & (\bm{R}_2)_{k,\ell} &\coloneqq \begin{cases} r_2, & k = \ell = 2j - 1, \; j = 1, \ldots, N_2, \\ 0, & \text{otherwise}, \end{cases} \\
    \bm{C}_{k,\ell} &\coloneqq \begin{cases} -1, & k = 2N_1 + 1, \; \ell = 1, \\ 0, & \text{otherwise}, \end{cases} & \Hamiltonian(\x) &= \sum_{i=1}^2 \sum_{j=1}^{N_i} \left( \frac{p_{i,j}^2}{2m_i} + \frac{K_i}{2} q_{i,j}^2 + \frac{\beta_i}{4} q_{i,j}^4 \right) + \frac{K_c}{2} q_c^2,
\end{align*}
where $\bm{e}_1 = (1,0,\ldots,0)^\top \in \mathbb{R}^{2N_1}$ and $\operatorname{tridiag}(a,b,c)$ denotes the tridiagonal matrix with the entries $a$, $b$, $c$ on the sub-, main-, and superdiagonal.
Throughout this section, we consider
\begin{equation*}
    \begin{aligned}
        m_1 &= 10^{-2}, & m_2 &= 1, & K_1 &= 25, & K_2 &= 1, & \beta_1 &= 10^3, & \beta_2 &= 4, & K_c &= 0.1, & r_1,r_2 &\geq 0.
    \end{aligned}
\end{equation*}
Following \Cref{remark:coupled_pHS}, we separate the dissipative part $\{3\}$ and split $\bm{J}$ into the slow subsystem $\{1\}$ and the fast subsystem $\{2\}$, i.e., $\x' = \bm{f}^{\{m\}}(\x) = \bm{A}^{\{m\}}\nabla\Hamiltonian(\x)$, with
$$ \bm{A}^{\{1\}} = \begin{pmatrix}
    \bm{0} & \bm{0} \\ \bm{0} & \bm{J}_2
\end{pmatrix}, \quad \bm{A}^{\{2\}} = \begin{pmatrix}
    \bm{J}_1 & \bm{C} \\ -\bm{C}^\top & \bm{0}
\end{pmatrix}, \quad \bm{A}^{\{3\}} = \begin{pmatrix}
    -\bm{R}_1 & \bm{0} \\ \bm{0} & -\bm{R}_2
\end{pmatrix}. $$
As $\Hamiltonian$ is quadratic in $\p$ and $\bm R$ is diagonal with $R_{k,k}=r_i$ for the momenta of chain $i$ and $0$ otherwise, $\x^\prime=-\bm R\nabla\Hamiltonian(\x)=\bm\Lambda\x$ with $\lambda_k=-r_i/m_i$ for the momenta of chain $i$ and $\lambda_k=0$ otherwise.
Hence, the flow of the dissipative subsystem $\{3\}$ can be computed exactly:
\begin{equation}\label{eq:R_flow_exact}
    \flow{\stepsize}{3}(\x_0) = \exp( \stepsize \bm{\Lambda} ) \x_0 .
\end{equation}
Numerical simulations are performed on the interval $[0,10]$ with initial values
\begin{equation}\label{eq:MSD_initial_values}
    \begin{aligned}
        q_{1,j}(0) &= \tfrac{1}{50} \sin\big( \tfrac{\pi j}{N_1 + 1} \big), &
        q_{2,j}(0) &= \tfrac{3}{10} \, \mathrm{e}^{-((j - 7)/2)^2} + \tfrac{1}{10} \sin\big( \tfrac{15 \pi j}{N_2 + 1} \big), \\
        p_{1,j}(0) &= 0, &
        p_{2,j}(0) &= -\tfrac{3}{10} \, \mathrm{e}^{-((j - 7)/2)^2}, \quad q_c(0)=0.
    \end{aligned}
\end{equation}

Recalling \Cref{sec:design_PHS}, an energy-consistent scheme is required to reproduce the correct qualitative energy behavior.
Conventional splitting methods of order $\convergenceorder > 2$ necessarily involve negative integrator coefficients \cite{BLANES200523}, thereby violating the energy consistency. 
We consider the splitting tree depicted in \Cref{fig:MSD_splitting_tree} and assign the commutator-based splitting method \eqref{eq:JR4} to the root node, to the node $\{1,2\}$ the triple-jump composition \eqref{eq:Yoshida_compressed} of the Strang splitting with the roles of $\{1\}$ and $\{2\}$ interchanged, the modified discrete gradient method \eqref{eq:DGM4} to the nodes $\{1\}$ and $\{2\}$, and the exact flow \eqref{eq:R_flow_exact} to node $\{3\}$.
We refer to the resulting hierarchical splitting method as \texttt{HSM4}. 
Since all methods assigned to nodes of the splitting tree are convergent of order four, the overall hierarchical splitting method has order $\convergenceorder=4$.
Moreover, the method is energy-consistent in the sense of \Cref{sec:design_PHS}: all (perturbed) subsystems are pH-ODEs with respect to $\Hamiltonian$, the dissipative subsystem is only integrated forward in time (the negative weight $\gamma_2$ of the triple-jump composition enters only the conservative subsystems), and the conservative subsystems are integrated with an energy-consistent discrete gradient method.
Consequently, it holds $\Hamiltonian(\x_{n+1}) \le \Hamiltonian(\x_n)$ for all step sizes $\stepsize > 0$.

As reference methods, we consider \texttt{BMs6p4}, the composition \eqref{eq:BMs6p4} with base scheme
\begin{equation*}
    \bm{\chi}_{\stepsize} = \flow{\stepsize}{3} \circ \big( \SplittingMethod{\stepsize/\mrfactor{}}{2} \big)^{\mrfactor{}} \circ \SplittingMethod{\stepsize}{1},
\end{equation*}
where the numerical flows of the subsystems $\{1\}$ and $\{2\}$ are computed by \eqref{eq:DGM4}.
The method \eqref{eq:BMs6p4} includes negative coefficients, such that the dissipative subsystem $\{3\}$ is evaluated backward in time, and is thus not energy-consistent.
The second reference method, \texttt{DGM4}, is the fourth-order energy-consistent modified discrete gradient method \eqref{eq:DGM4} applied to the monolithic pH-ODE \eqref{eq:coupled_MSD}, i.e., without splitting and without multiple time stepping. 
In all discrete gradient methods, the average vector field discrete gradient is evaluated exactly by the two-point Gauss--Legendre rule, as $\nabla\Hamiltonian$ is a cubic polynomial.
All nonlinear equations are solved by Newton's method with the analytic Jacobian and dense LU factorization, using a termination tolerance of $10^{-12}$.
The maximum number of iterations is limited to $50$.

To illustrate the importance of energy consistency, we first consider a stiffly damped configuration with $(r_1,r_2) = (20,0)$.
Due to the strong dissipation in the fast chain, the system offers no multirate potential to exploit.
We therefore set $\mrfactor{} = 1$ for \texttt{BMs6p4} and \texttt{HSM4}.
As $\Hamiltonian$ is non-increasing along the exact solution, a numerical solution should satisfy $\max_n [\Hamiltonian(\x_{n+1}) - \Hamiltonian(\x_n)]/\stepsize \leq 0$.
\Cref{tab:MSD_energy} shows that \texttt{HSM4} and \texttt{DGM4} keep the Hamiltonian non-increasing for all step sizes considered, whereas \texttt{BMs6p4} violates the dissipation inequality for $\stepsize = 5/512$ and, for larger step sizes, the Newton iteration fails to converge.
This is caused by the sub-steps $\x \gets \exp(c \stepsize \bm{\Lambda}) \x$ with the negative coefficient $c = \alpha_3 = \beta_4 \approx -0.41$ of the composition, which amplify the stiffly damped momenta and thereby increase the Hamiltonian.
Hence, \texttt{BMs6p4} is only conditionally stable.

\begin{figure}[htbp!]
    \centering
    \begin{minipage}[c]{.3\textwidth}
        \centering
        \begin{tikzpicture}[level distance=0.8cm,
              level 1/.style={sibling distance=2cm},
              level 2/.style={sibling distance=2cm},
              level 3/.style={sibling distance=2cm}]
              \node {$\{1,2,3\}$}
                child {
                node {$\{1,2\}$}
                  child {node {$\{1\}$}
                  }
                  child {node {$\{2\}$}
                  edge from parent node[sloped,above,draw=none] {\tiny~\;$\mrfactor{}$}
                  }
                }
                child {node {$\{3\}$
                }
                };
        \end{tikzpicture}
    \end{minipage}
    \hfill
    \begin{minipage}[c]{.65\textwidth}
        \centering
        \small
        \setlength{\tabcolsep}{2.5pt}
        \begin{tabular}{@{}c|c|ccc@{}}
            \toprule
            $\stepsize$ & exact & \texttt{BMs6p4} & \texttt{HSM4} & \texttt{DGM4} \\
            \midrule
    $5/2$ & $-1.4\cdot 10^{-3}$ & failed & 0 & $-4.7\cdot 10^{-7}$ \\
    $5/32$ & $-9.8\cdot 10^{-4}$ & failed & $-1.8\cdot 10^{-12}$ & $-1.7\cdot 10^{-4}$ \\
    $5/256$ & $-9.6\cdot 10^{-4}$ & failed & $-8.2\cdot 10^{-6}$ & $-2.4\cdot 10^{-10}$ \\
    $5/512$ & $-9.6\cdot 10^{-4}$ & $1.3\cdot 10^{2}$ & $-6.6\cdot 10^{-10}$ & $-9.3\cdot 10^{-9}$ \\
    $5/1024$ & $-9.5\cdot 10^{-4}$ & $-1.0\cdot 10^{-3}$ & $-7.7\cdot 10^{-4}$ & $-1.5\cdot 10^{-7}$ \\
    $5/2048$ & $-9.5\cdot 10^{-4}$ & $-9.6\cdot 10^{-4}$ & $-1.1\cdot 10^{-3}$ & $-1.2\cdot 10^{-4}$ \\
    $5/8192$ & $-9.5\cdot 10^{-4}$ & $-9.5\cdot 10^{-4}$ & $-9.5\cdot 10^{-4}$ & $-9.6\cdot 10^{-4}$ \\
            \bottomrule
        \end{tabular}
    \end{minipage}

    \begin{minipage}[t]{.3\textwidth}
        \captionof{figure}{Splitting tree of the hierarchical splitting method for the coupled MSD chains \eqref{eq:coupled_MSD}.}
        \label{fig:MSD_splitting_tree}
    \end{minipage}
    \hfill
    \begin{minipage}[t]{.65\textwidth}
        \captionof{table}{Largest discrete energy change rate $\max_n[\Hamiltonian(\x_{n+1})-\Hamiltonian(\x_n)]/\stepsize$ for $(r_1,r_2)=(20,0)$, $\mrfactor{}=1$ (entries of modulus below $10^{-15}$ shown as $0$).
        ``exact'': reference solution on the same grid; ``failed'': Newton iteration did not converge.}
        \label{tab:MSD_energy}
    \end{minipage}
\end{figure}

\begin{figure}[t]
    \centering
    \includegraphics{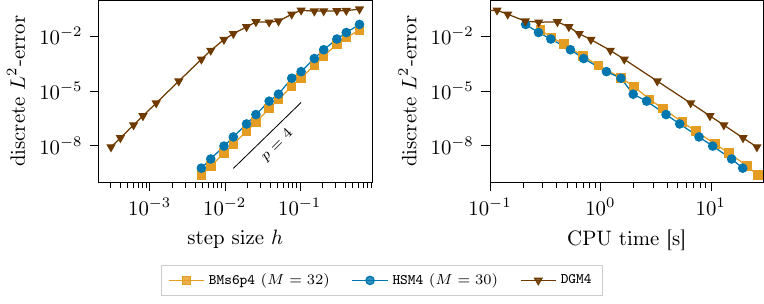}
    \caption{Numerical results for coupled MSD chains with $(r_1,r_2) = (10^{-4},10^{-2})$. 
    Left: discrete $L^2$-error vs.~step size $\stepsize$. 
    Right: discrete $L^2$-error vs.~CPU time (in seconds).}
    \label{fig:MSD_results}
\end{figure}

For testing the accuracy and computational efficiency of the methods, we consider the weakly damped configuration $(r_1,r_2) = (10^{-4},10^{-2})$.
For this choice, the first chain is light and stiff with linearized frequency $\sqrt{K_1/m_1} = 50$, whereas the second chain is heavy and soft with linearized frequency $\sqrt{K_2/m_2} = 1$.
This multirate potential can be exploited in \texttt{BMs6p4} and \texttt{HSM4} via the multirate factor $\mrfactor{}$, where both variants use the reweighting \eqref{eq:reweighted_mr-factor}.
We choose $\mrfactor{} = 32$ for \texttt{BMs6p4} and $\mrfactor{} = 30$ for \texttt{HSM4}, obtained by manual tuning so as to achieve the best performance in the work-precision diagram below.
The error of a numerical solution $\x_0,\ldots,\x_{N_t}$ with $\stepsize = 10/N_t$ is measured in the discrete $L^2$-norm $\big( \stepsize \sum_{n=0}^{N_t} \lVert \x_n - \x(t_n) \rVert_2^2 \big)^{1/2}$.
\Cref{fig:MSD_results} (left) confirms that all three methods achieve the expected convergence order $\convergenceorder = 4$.
Note that \texttt{DGM4} attains its order only for smaller step sizes, since its global step size has to resolve the fast chain.
The work-precision diagram in \Cref{fig:MSD_results} (right) shows that exploiting the multirate potential via multiple time stepping yields a more efficient computational process: \texttt{HSM4} requires about $10\,\%$ less CPU time than \texttt{BMs6p4} and roughly $50\,\%$ of the CPU time of \texttt{DGM4}.
The efficiency gain over \texttt{BMs6p4} is thus moderate and depends on the manually tuned multirate factors.

Since \texttt{BMs6p4} is itself a hierarchical splitting method by \Cref{lemma:composition_as_HSM}, the comparison contrasts two methods of order four that differ in the splitting tree and the methods assigned to the nodes.
The example satisfies all three conditions of \Cref{remark:superiority}: the block structure provides multirate potential, which \texttt{DGM4} cannot exploit, commutator-based methods are applicable, and assigning an energy-consistent method to a node retains the dissipation inequality, which the composition violates.

\section{Conclusions}
In this paper, we introduced a novel class of hierarchical splitting methods for $\Npartitions$-split differential equations.
These methods apply splitting methods for two-split systems in a hierarchical manner, organized by a splitting tree.
Since assigning numerical integration schemes of order $\convergenceorder$ to all nodes of a splitting tree suffices to obtain a hierarchical splitting method of order $\convergenceorder$, the framework simplifies the construction of splitting methods for $\Npartitions$-split systems that are not based on composition techniques.
We derived explicit formulas for the leading-order error terms, investigated self-adjointness, and showed that compositions of hierarchical splitting methods are again hierarchical splitting methods, such that state-of-the-art splitting methods for $\Npartitions$-split systems are contained in the proposed class as special cases.

We also introduced an extension to multirate hierarchical splitting methods.
The resulting multiple time stepping schemes retain the structure-preserving properties and inherit the order theory of their singlerate counterparts.
We augmented the multirate extension by a reweighting technique for the multirate factors to reduce computational overhead, and introduced the concept of computational order, providing conditions on the multirate factors under which superior accuracy is expected over a restricted step size range.
At the same time, we have shown that the reweighting can lead to an order reduction for composition techniques, such that constant multirate factors must be employed to retain the convergence order.
Based on this framework, we discussed principles for designing splitting trees, with particular attention to separable Hamiltonian systems and port-Hamiltonian systems.
For separable Hamiltonian systems, embedding nested integration schemes into the hierarchical splitting approach yields a refined error analysis, reduces computational overhead via the reweighting, and makes force-gradient integrators available for the inner nodes of a splitting tree.
For port-Hamiltonian systems, the hierarchical splitting approach allows for the construction of energy-consistent splitting methods of order four that simultaneously exploit multiple time stepping.

The numerical results confirm the theoretical findings and quantify the resulting gains.
For the Fermi--Pasta--Ulam problem, the reweighting of the multirate factors reduces the computational cost of the hierarchical splitting methods while preserving their convergence order, whereas it reduces the order of composition techniques to two.
Among the methods of order four, the hierarchical splitting methods are computationally more efficient and more robust for increasing step sizes than the composition schemes, and embedding the force-gradient approach leads to the most efficient method in this benchmark.
For the coupled mass-spring-damper chains, the hierarchical splitting method combines the order and the energy consistency of a discrete gradient method with multiple time stepping, whereas the composition scheme violates the dissipation inequality and is only conditionally stable.

The hierarchical splitting framework inherits accuracy, efficiency, and structure preservation from the schemes assigned to its nodes, and therefore benefits directly from any advance in splitting methods for two-split systems and in numerical integration schemes for the individual subsystems.

The hierarchical splitting approach does not resolve how to decompose a given vector field into $\Npartitions$ partitions, nor how to order the operators within a splitting method in general. 
Both remain open challenges in the numerical time integration of partitioned systems.

We expect the structural benefits of hierarchical splitting to become even more relevant for large-scale systems demanding larger multirate factors, such as multi-physics systems that are potentially modeled as port-Hamiltonian systems \cite{VanDerSchaft2014PHS}.
A systematic treatment of the block structure of the system matrices of such coupled systems is beyond the scope of the present work, and we will thus, as a next step, investigate hierarchical splitting approaches for (coupled) port-Hamiltonian systems in more detail.

While the embedding of nested integration schemes into the hierarchical splitting approach leads to a natural extension of the order theory, a corresponding extension of the numerical stability analysis available for two-split systems \cite{blanes2008stability,SCHAFERS2026110034} is still missing.
Deriving step size restrictions for the individual nodes of a splitting tree and their interplay with the multirate factors is therefore a natural subject of a separate study.

Furthermore, the hierarchical splitting approach presents a promising strategy for the derivation of splitting methods for $\Npartitions$-split systems beyond the ODE case, such as partial differential equations (PDEs)~\cite{blanes2024splitting}, differential-algebraic equations (DAEs)~\cite{BARTEL202528}, and stochastic differential equations (SDEs)~\cite{Foster2024splitting}.
However, depending on the specific class of equations, a more comprehensive investigation of the hierarchical splitting approach is necessary.
For instance, it remains an open question to what extent splitting methods with complex coefficients~\cite{castella_splitting_2009,HansenOstermann2009,blanes2025splitting_complex} can be integrated into the hierarchical splitting approach to obtain methods whose coefficients have positive real part.

\backmatter

\bmhead{Acknowledgements}
This work is supported by the Deutsche Forschungsgemeinschaft (DFG, German Research Foundation) research unit FOR5269 \enquote{Future methods for studying confined gluons in QCD}.
Michael Günther acknowledges funding by the DFG -- Project-ID 531152215 -- CRC 1701.
The authors thank Andreas Bartel for helpful comments on an earlier version of this manuscript.

\begin{appendices}

\section{Numerical Methods}
\texttt{BMs6p4}. Given a base scheme $\bm{\chi}_{\stepsize}$, a fourth-order composition method of the form \eqref{eq:composition_based_on_LT} is obtained for $s=6$ and \cite{Blanes_Moan2002}
\begin{equation}\label{eq:BMs6p4}
    \begin{aligned}
        \beta_1 &= \alpha_6 = 0.162314550766870, & 
        \alpha_1 &= \beta_6 = 0.0829844064174052, \\
        \beta_2 &= \alpha_5 = 0.370877414979582, & 
        \alpha_2 &= \beta_5 = 0.233995250731498, \\
        \beta_3 &= \alpha_4 = 0.059762097006575, & 
        \alpha_3 &= \beta_4 = -0.409933719901930.
    \end{aligned}
\end{equation}

\noindent\texttt{BMs10p6}. Given a base method $\bm{\chi}_{\stepsize}$, a sixth-order composition method of the form \eqref{eq:composition_based_on_LT} is obtained for $s=10$ and \cite{Blanes_Moan2002} 
\begin{equation}\label{eq:BMs10p6}
    \begin{aligned}
        \alpha_1 &= \beta_{10} = 0.05026276440039221, & 
        \alpha_2 &= \beta_9 = 0.31496061692769417, \\
        \alpha_3 &= \beta_8 = 0.49242637248987586, & 
        \alpha_4 &= \beta_7 = 0.23706391397812188, \\
        \alpha_5 &= \beta_6 = 0.34635818985072686, &
        \alpha_6 &= \beta_5 = -0.36276277925434486, \\
        \alpha_7 &= \beta_4 = 0.19560248860005314, & 
        \alpha_8 &= \beta_3 = -0.42511876779769087, \\
        \alpha_9 &= \beta_2 = -0.44734648269547816, &
        \alpha_{10} &= \beta_1 = 0.0985536835006498.
    \end{aligned}
\end{equation}

\noindent\texttt{JR4}. Given a pH-ODE of the form 
\begin{equation*}
\begin{aligned}
    \begin{pmatrix}
        \p^\prime \\ \q^\prime
    \end{pmatrix} &= (\bm{J} - \bm{R}) \begin{pmatrix}
        \nabla_{\p} \Hamiltonian(\p,\q) \\ \nabla_{\q} \Hamiltonian(\p,\q)
    \end{pmatrix}, & \bm{R} &= \begin{pmatrix}
        \tilde{\bm{R}} & \bm{0} \\ \bm{0} & \bm{0}
    \end{pmatrix}, & \Hamiltonian(\p,\q) &= \frac{1}{2} \p^\top \bm{M}^{-1} \p + \Potential(\q), 
\end{aligned}
\end{equation*}
with $\bm{M} = \bm{M}^\top$, an energy-consistent splitting method of order four is given by \cite{Monch2026splitting}
\begin{equation}\label{eq:JR4}
    \SplittingMethod{\stepsize}{} = \flow{\stepsize/6}{\bm{R}} \circ \tilde{\bm{\varphi}}_{\stepsize/2}^{\{\bm{J}\}} \circ \flow{2\stepsize/3}{\bm{R}} \circ \tilde{\bm{\varphi}}_{\stepsize/2}^{\{\bm{J}\}} \circ \flow{\stepsize/6}{\bm{R}}, 
\end{equation}
where $\flow{\stepsize}{\bm{R}}$ and $\tilde{\bm{\varphi}}_{\stepsize}^{\{\bm{J}\}}$ denote the exact flows of the (perturbed) subsystems
\begin{equation*}
        \begin{pmatrix} \p^\prime \\ \q^\prime \end{pmatrix} \!=\! \begin{pmatrix}
            - \tilde{\bm{R}} \bm{M}^{-1} \p \\ \bm{0}            
        \end{pmatrix} \, \text{and} \,     
        \begin{pmatrix}
            \p^\prime \\ \q^\prime
        \end{pmatrix} \!=\! \big(\!\bm{J} - \tfrac{\stepsize^2}{72} \big[ \bm{N} , [ \bm{N} , \bm{J}] \big] \!\big) \!\begin{pmatrix}
            \nabla_{\p} \Hamiltonian(\p,\q) \\ \nabla_{\q} \Hamiltonian(\p,\q)
        \end{pmatrix}, \; \bm{N} \!\coloneqq\! \begin{pmatrix}
    -\tilde{\bm{R}} \bm{M}^{-1} & \bm{0} \\ \bm{0} & \bm{0}
\end{pmatrix}.
\end{equation*}

\noindent\texttt{DGM4}. For closed pH-ODEs $\x^\prime = (\bm{J} - \bm{R}) \nabla \Hamiltonian(\x)$, the modified discrete gradient method \cite{Lautwein2026DGM}
\begin{equation}\label{eq:DGM4}
    \begin{aligned}
        \frac{\x_1 - \x_0}{\stepsize} &= \tilde{\bm{A}}(\x_0 , \x_1 , \stepsize) \overline{\nabla}_{\mathrm{AVF}} \Hamiltonian(\x_0,\x_1), \\
        \tilde{\bm{A}}(\x_0,\x_1,\stepsize) &= (\bm{J} - \bm{R}) - \tfrac{\stepsize^2}{12} (\bm{J}-\bm{R}) \bm{H} (\bm{J}-\bm{R}) \bm{H} (\bm{J}-\bm{R}) - \stepsize^4 \bm{D}^\top \bm{R} \bm{D}, \\
        \bm{H} &= \nabla^2 \Hamiltonian \big( \tfrac{\x_0 + \x_1}{2} \big), \quad
        \bm{D} = -\tfrac{1}{12} \bm{H} \bm{J} \bm{H} \bm{J} -\tfrac{1}{24}\bm{H}\bm{R}\bm{H} \bm{R} , \\
        \overline{\nabla}_{\mathrm{AVF}} \Hamiltonian(\x_0,\x_1) &= \int_{0}^{1} \nabla \Hamiltonian((1-\xi)\x_0+\xi \x_1) \, \mathrm{d}\xi ,
    \end{aligned}
\end{equation}
is energy-consistent and has order $\convergenceorder = 4$. For conservative systems ($\bm{R} \equiv \bm{0}$), the correction matrix simplifies to 
\begin{equation}\label{eq:DGM4_J}
    \begin{aligned}
        \tilde{\bm{A}}(\x_0,\x_1,\stepsize) = \bm{J} - \tfrac{\stepsize^2}{12} \bm{J} \bm{H} \bm{J} \bm{H} \bm{J}.
    \end{aligned}
\end{equation}
For separable Hamiltonian systems $\Hamiltonian(\p,\q) \!=\! \Kinetic(\p) + \Potential(\q)$, we consider splitting methods 
\begin{equation}\label{eq:splitting_method_Hamiltonian}
    \SplittingMethod{\stepsize}{} = \flow{b_s \stepsize}{\Potential} \circ \flow{a_s \stepsize}{\Kinetic} \circ \ldots \circ \flow{b_1 \stepsize}{\Potential} \circ \flow{a_1 \stepsize}{\Kinetic},
\end{equation}
with $\flow{\stepsize}{\Kinetic}$ and $\flow{\stepsize}{\Potential}$ denoting the exact flow of the kinetic and potential part, respectively.

\noindent\texttt{OMF4}. A particularly efficient six-stage method of order $\convergenceorder = 4$ is defined by \cite{omelyan2003symplectic}
\begin{equation}\label{eq:OMF4}
\begin{aligned}
    a_1 &= 0, & b_1 &= b_6 = 0.083983152628767, \\
    a_2 &= a_6 = \phantom{-}0.253978510841060, & b_2 &= b_5 = 0.682236533571909, \\
    a_3 &= a_5 = -0.032302867652700, & b_3 &= b_4 = 0.5-(b_1+b_2), \\
    a_4 &= 1-2(a_2+a_3).
\end{aligned}
\end{equation}

\noindent\texttt{OMF6}. A particularly efficient eight-stage method of order $\convergenceorder = 6$ is obtained for \cite{omelyan2003symplectic}
\begin{equation}\label{eq:OMF6}
    \begin{aligned}
        a_1 &= 0, & b_1 &= b_8 = \phantom{-}0.0833333333333333, \\
        a_2 &= a_8 = \phantom{-}0.2465881872786138, & b_2 &= b_7 = \phantom{-}0.3977675859548440, \\
        a_3 &= a_7 = \phantom{-}0.6047073875057809, & b_3 &= b_6 = -0.0393336931446257, \\
        a_4 &= a_6 = -0.4009869039788007, & b_4 &= b_5 = 0.5 - (b_1 + b_2 + b_3), \\
        a_5 &= 1 - 2(a_1 + a_2 + a_3 + a_4). & & 
    \end{aligned}
\end{equation}

For separable Hamiltonian systems with quadratic kinetic energy \eqref{eq:separable_H}, we can replace the momentum updates in \eqref{eq:splitting_method_Hamiltonian} by force-gradient updates 
$$ \flow{b_i \stepsize, c_i \stepsize^3}{\Potential,\mathcal{C}} = \exp \left( b_i \stepsize \LieDerivative{\Potential} + c_i \stepsize^3 \LieDerivative{\mathcal{C}}  \right) \mathrm{Id}, $$
with the Lie derivatives 
$$ \LieDerivative{\Potential} = - \nabla \Potential(\q) \frac{\partial}{\partial \p} \quad \text{and} \quad \LieDerivative{\mathcal{C}} = \nabla \Potential(\q)^\top \bm{M}^{-1} \nabla^2\Potential(\q) \frac{\partial}{\partial \p},$$
thereby obtaining force-gradient integrators 
\begin{equation}\label{eq:FGI}
    \SplittingMethod{\stepsize}{} = \flow{b_s \stepsize, c_s \stepsize^3}{\Potential,\mathcal{C}} \circ \flow{a_s \stepsize}{\Kinetic} \circ \ldots \circ \flow{b_1 \stepsize, c_1 \stepsize^3}{\Potential,\mathcal{C}} \circ \flow{a_1 \stepsize}{\Kinetic},
\end{equation}

\noindent\texttt{ABACABACABA}. In the classification of force-gradient integrators conducted in \cite{omelyan2003symplectic}, the most efficient fourth-order force-gradient integrator is obtained for $s=6$ and coefficients
\begin{equation}\label{eq:ABACABACABA}
    \begin{aligned}
        a_1 &= a_6 = 0.0641910886681624, & b_1 &= b_5 = 0.1518179640276466, \\
        a_2 &= a_5 = 0.1919807940455741, & b_2 &= b_4 = 0.2158369476787619, \\
        a_3 &= a_4 = 0.5 - (a_1 + a_2),  & b_3 &= 1 - 2(b_1 + b_2), \\
        c_2 &= c_4 = 0.0009628905212025, & b_6 &= c_1 = c_3 = c_5 = c_6 = 0.
    \end{aligned}
\end{equation}

\end{appendices}

{
\setlength{\bibsep}{0pt}
\renewcommand*{\bibfont}{\scriptsize}

}

\end{document}